\documentclass[reqno, 12pt]{amsart}
\usepackage{amsfonts}
\usepackage{amssymb,amsmath}
\usepackage{amsthm}
\usepackage{upref}
\usepackage{enumerate}
\usepackage{bbm}

\makeatletter
\def\LaTeX{\leavevmode L\raise.42ex
    \hbox{\kern-.3em\size{\sf@size}{0pt}\selectfont A}\kern-.15em\TeX}
 \makeatother

 \DeclareMathOperator*{\slim}{s-lim}
\numberwithin{equation}{section}

\newtheorem{lemma}{Lemma}[section]
\newtheorem{theorem}[lemma]{Theorem} 
\newtheorem{corollary}[lemma]{Corollary}
\newtheorem{proposition}[lemma]{Proposition}

\theoremstyle{definition}

\newtheorem{definition}[lemma]{Definition}
\newtheorem{example}[lemma]{Example}

\newtheorem{assumption}[lemma]{Assumption}

\newtheorem{remark}[lemma]{Remark}

  \newcommand{\e}{\eqref}
\newcommand{\ri}{\rightarrow}
\newcommand{\q}{\quad}

\newcommand{\wt}{\widetilde}
\newcommand{\wh}{\widehat}
\newcommand{\la}{\langle}
\newcommand{\ra}{\rangle}

\newcommand{\ov}{\overline}
 \renewcommand{\d}{\delta}

   \newcommand{\sgn}{\operatorname{sgn}}
  \newcommand{\ran}{\operatorname{Ran}}

\renewcommand\Im{\operatorname{Im}}

\newenvironment{pf}{\begin{proof}}{\end{proof}}

\def\qqq{\mathrel{\subset\mkern-15mu\lower.38ex\hbox{${\scriptscriptstyle\rightarrow}$}}}
\let\cal\mathcal

\let\Bbb\mathbb

 \DeclareMathOperator{\Ran}{Ran}
  
  \DeclareMathOperator{\spec}{spec}
  
\begin{document}

\title {Spectral and  scattering theory for  differential and Hankel  operators}

%\title {Spectral and  scattering theory for  quasi-Carleman operators}
%\title {Spectral   theory of   Carleman type Hankel operators}
\author{ D. R. Yafaev}
\address{ IRMAR, Universit\'{e} de Rennes I\\ Campus de
  Beaulieu, 35042 Rennes Cedex, FRANCE}
\email{yafaev@univ-rennes1.fr}

%  \date{\today}

\begin{abstract}
We   consider a class of Hankel operators $H$ realized
  in the space $L^2 ({\Bbb R}_{+}) $ as
integral operators  with  kernels $h(t+s)$ where $h(t)=P (\ln t ) t ^{-1}$ and $P(X)=   X^n+p_{n-1}  X^{n-1}+\cdots$ is an arbitrary real polynomial   of degree $n$. This class contains the classical Carleman operator when $n =0$.  We show that a Hankel operator $H$ in this class can be reduced by an {\it explicit}  unitary transformation
 (essentially by the Mellin transform)  
 to a differential operator $A = v Q(D) v$ in the space $L^2 ({\Bbb R}) $. Here $Q(X)=   X^n+ q_{n-1}  X^{n-1}+\cdots$ 
 is a polynomial determined by $P(X)$ and $v(\xi)=\pi^{1/2} (\cosh(\pi\xi))^{-1/2} $ is the universal function. Then the
 operator $A = v Q(D) v$  reduces by the generalized Liouville transform to the standard differential operator 
$B =    D^n+ b_{n-1}  (x)D^{n-1}+\cdots+ b_{0} (x)$  with the coefficients $b_{m}(x)$, $m=0,\ldots, n-1$,   decaying  sufficiently rapidly as $|x|\to \infty$. This allows us to use the results of spectral   theory of differential operators for the study of spectral properties of generalized Carleman operators. In particular, we show that the absolutely continuous spectrum of $H$ is simple and coincides with $\Bbb R$ if $n$ is odd, and it has multiplicity $2$ and coincides with $[0,\infty)$ if $n\geq 2$ is even. The singular continuous spectrum of $H$ is empty, and its eigenvalues may accumulate to the point $0$ only. As a by-product of our considerations, we develop spectral theory of a new class of {\it degenerate} differential operators  $A = v Q(D) v$  where $Q(X)$ is an  arbitrary real polynomial  and $v(\xi)$ is a sufficiently arbitrary real function decaying at infinity.
   \end{abstract}

\maketitle

% \thispagestyle{empty}

%************************************************************
\section{Introduction}  
%***********************************************************

{\bf 1.1.}
Hankel operators  can be defined by the formula
\begin{equation}
(H u)(t) = \int_{0}^\infty h(t+s) u(s)ds 
\label{eq:H1}\end{equation}
in the space $L^2 ({\Bbb R}_{+}) $; thus integral kernels $h$ of Hankel operators  depend  on the sum of variables only. 
We refer to the books \cite{NK, Pe} for basic information on Hankel operators. Of course $H$ is symmetric if $  h(t)=\ov{h(t)}$.  

The spectra of bounded self-adjoint   Hankel operators were characterized in the fundamental paper \cite{MPT}. 
 It was shown in \cite{MPT} that the spectrum of a Hankel  operator   contains the point zero, and if zero is an eigenvalue, then necessary it has infinite multiplicity. Moreover,   the spectral multiplicities of the points $\lambda$ and $-\lambda$ cannot differ by more than $2$ and they cannot differ by more than  $1$ on the singular spectrum.
 Conversely, if the spectral  measure and the multiplicity function of  a self-adjoint operator  $H$ possess these properties, then $H$ is unitarily equivalent to a Hankel operator.

 Since this  result applies to all self-adjoint   Hankel operators, it does not allow one to find spectral properties of specific
 classes of Hankel operators.  The  cases where Hankel operators can be explicitly diagonalized are very scarce.     The simplest and most important kernel $h(t)=t^{-1}$  was considered by T.~Carleman in \cite{Ca}. The eigenfunctions of the continuous spectrum $\theta (t , k)$, $k\in {\Bbb R}$, of this operator are given by the formula $\theta(t , k) = t^{-1/2 + ik}$. They satisfy the equation $H \theta(k)=\lambda(k) \theta(k)$ with the dispersion relation $\lambda(k)= \pi \big(\cosh (\pi k)\big)^{-1}$. Thus the spectrum of $H$ is absolutely continuous, it has multiplicity $2$ and coincides with the interval $[0,\pi]$. It was pointed out by J.~S.~Howland in   \cite{Howland} that there is a somewhat mysterious affinity between Hankel and differential operators. In terms of this analogy, the Carleman operator plays the role of the operator $ D^2$ in the space $L^2({\Bbb R})$.
 
 The results on  the Carleman operator can be extended to  more complicated   kernels. We note the classical papers
 \cite{Me} by  F.~Mehler who considered the kernel $h(t)=(t+1)^{-1}$ and \cite{Ma,Ro} by W.~Magnus   and   M.~Rosenblum who considered the kernel $h(t)=t^{-1}e^{-t}$ (M.~Rosenblum   considered also more general kernels with the same singularity at $t=0$). The corresponding Hankel operators $H$ were diagonalized in terms of the Legendre and Whittaker functions, respectively. The spectrum of these operators is absolutely continuous, simple,  and it coincides with the interval $[0,\pi]$. These results can be deduced from the fact that Hankel operators with such kernels commute with  appropriate differential operators of second order; see \cite{Y1} where some new examples are also considered. 
  
  Another generalization of the Carleman operator is given by the kernel $h(t)=t^{-q}$, $q>0$. It was shown in
  \cite{Yf1a} that, for $q\neq 1$,  the spectrum of the corresponding Hankel operator $H$  is absolutely continuous, it has constant multiplicity (presumably, it is simple) and coincides with the half-axis $[0, \infty)$. The approach of
  \cite{Yf1a} relies only on the invariance property of $H$ with respect to the group of the dilations $f(t)\mapsto \sqrt{\rho}f(\rho t)$, $\rho>0$. So it does not give any information on the structure of   eigenfunctions of the operator $H$.
  
  We also  note  the paper   \cite{K-M} by V.~Kostrykin and K.~Makarov    where it was shown that the Hankel operator $H$ with  kernel $h(t)=t^{-1} \sin t$ can be reduced by an appropriate unitary transformation to the orthogonal sum of two operators considered in  \cite{Ro}. The spectrum of $H$ is simple, absolutely continuous and coincides with $[-\pi/2,\pi/2]$.

\medskip
    
     {\bf 1.2.}
     Our goal      here is to study spectral properties of generalized Carleman operators with kernels  
 \begin{equation}
h(t)= P (\ln t)t^{-1}
\label{eq:LOG}\end{equation} 
where 
 \begin{equation}
P (X)=\sum_{m=0}^n p_mX^m
\label{eq:LOG4}\end{equation}
 is an arbitrary real polynomial. Without loss of generality we suppose that $ p_n =1$. Obviously, kernels \e{eq:LOG} have two singular points $t=\infty$ and $t=0$.
For $n\geq 1$, such  Hankel operators   are unbounded but are well defined as self-adjoint operators.
A large part of our results on generalized Carleman operators can be summarized by the following assertion.     Below we denote by $\la x\ra$ the operator of multiplication by the function $(1+x^2)^{1/2}$.

  \begin{theorem}\label{SpThHa}
Let $H $ be the self-adjoint Hankel   operator defined by formula \e{eq:H1} where $h(t)$ is function  \e{eq:LOG}  and $P(X)$ is a  real polynomial  \e{eq:LOG4} of degree $n\geq 1$. 
Then
\begin{enumerate}[\rm(i)]
\item
The spectrum of the operator $H$ is absolutely continuous except eigenvalues that  may accumulate to zero and  infinity only.   
 \item
  The limiting absorption principle holds, that is, for any $\sigma>1/2$, the operator-valued  function 
      \begin{equation}
  \la \ln t\ra^{-\sigma} ( H-z)^{-1} \la \ln t\ra^{-\sigma}, \q \Im z \neq 0,
 \label{eq:LAPH}\end{equation}
 is H\"older continuous with exponent $\alpha<\sigma-1/2$ $($and $\alpha<1)$ up to the real axis,  except the eigenvalues of the operator $H$ and the point zero.
   \item
The 
 absolutely continuous  spectrum of the operator $H$ covers $\Bbb R$ and is simple for $n$ odd. It coincides with $[0,\infty)$ and has multiplicity $2$  for $n$ even. 
 \item
 If $n$ is odd, then the multiplicities of eigenvalues  of the operator $H$  are bounded by $(n-1)/2$. If $n$ is even, then the multiplicities of positive eigenvalues are bounded by $n/2-1$, and the multiplicities of negative eigenvalues are bounded by $n/2$. 
 % The point zero is not an eigenvalue.
\end{enumerate}
    \end{theorem}
    
    Clearly, this assertion is similar in spirit to the corresponding results for differential operators of Schr\"odinger type (see Theorem~\ref{SpTh}, below). The difference, however, is that the weight  is $  \la \ln t\ra^{-\sigma} $ in  \e{eq:LAPH} while it is $  \la   x\ra^{-\sigma} $ (also with $\sigma>1/2$) for the resolvents of differential operators.
    Thus the power scale for differential operators corresponds to the logarithmic scale for Hankel operators.
  Note that for perturbations of the classical Carleman operator (when $h(t)$ is close to $t^{-1}$) the assertions (i)  and (ii) of Theorem~\ref{SpThHa} were obtained in \cite{Y2}; for such operators the  absolutely continuous  spectrum coincides with $[0, \pi)$ and has multiplicity $2$.
     
     %   We also observe that the weight $  \la \ln t\ra^{-\sigma}$, $\sigma>1/2$,  in \e{eq:LAPH} is the same  as
%  in the formulation of the limiting absorption principle for the classical Carleman operator. 
    
    The proof of the  limiting absorption principle requires a study of eigenfunctions $\theta (t,k)$, $k\in \Bbb R$, of the continuous spectrum of the operator $H$. They satisfy the equation $H \theta (k)= \lambda \theta (k)$ where $\lambda=k^n$, and their precise definition (see Section~5) is similar to that for differential operators.
      We need a uniform (away from the eigenvalues  of $H$ and from the point $\lambda=0$) estimate
     \[
    | \theta(t ,k)|\leq C t^{-1/2}
    \]
     and a similar estimate on  differences $\theta(t ,k')-\theta(t ,k)$. 
     
     Actually, for the proof of these estimates we establish a much stronger result
       finding the asymptotic behavior of   $\theta(t,k)$ as $t\to \infty$ and as $t\to 0$. For the statement of the required result, we need to introduce the phase function $\omega (t,k)$, $t\in\Bbb R_{+}$. Its precise definition will be given in Section~6 by formulas \e{eq:STP2q} -- \e{eq:BFGt}. Here we note that $\omega (t,k)$ is determined by the coefficient $p_{n-1}$ only and
       \begin{equation} 
    \omega (t,k)= -  n \pi^{-1}  \ln t \big(\ln \big| \frac{\ln t}{k} \big|  + O(1)\big)
   \label{eq:omeg}\end{equation}
       as $t\to\infty$ and $t\to 0$. We prove  that, for odd $n$,
       \begin{equation}  
\theta(t,k) = \sqrt{\frac{n}{\pi | k|}}t^{-1/2} \big(s(\lambda)  e^{  i    \omega (t,k)}+  e^{ - i    \omega (t,k)}+o(1)\big)  
\label{eq:Sodd2H}\end{equation}
as $t\to\infty$ if $k>0$ and as $t\to 0$ if $k<0$. If  $t\to \infty$ but $k<0$ or  $t\to 0$ but $k>0$, then
     \begin{equation}  
\theta(t,k) =  o(t^{-1/2}) . 
\label{eq:So-}\end{equation}
  For even $n$, we prove that
 \begin{equation}\left\{\begin{array}{lcl} 
\theta (t,k)= \sqrt{\frac{n}{\pi k}} t^{-1/2} \big( s_{11}(\lambda) e^{  i     \omega (t,k)} + e^{ - i    \omega (t,k)}+ o(1) \big),\quad t\ri \infty,
\\
\theta (t,k)= \sqrt{\frac{n}{\pi k}} t^{-1/2} \big(     s_{21}(\lambda) e^{-  i    \omega (t,k)}+ o(1) \big),\quad t\ri 0.
\end{array}\right.
\label{eq:S1H}\end{equation}
if $k>0$, and
\begin{equation}\left\{\begin{array}{lcl} 
\theta(t,k)= \sqrt{\frac{n}{\pi | k|}} t^{-1/2} \big( s_{12}(\lambda) e^{  i    \omega (t,k)}+ o(1) \big),\quad t\ri \infty,
\\
\theta(t,k)= \sqrt{\frac{n}{\pi |k|}} t^{-1/2} \big(   e^{  i    \omega (t,k)} +  s_{22}(\lambda)  e^{-  i    \omega (t,k)}+ o(1) \big),\quad t\ri 0.
\end{array}\right.
\label{eq:S2H}\end{equation}
if $k<0$.  Here  $s(\lambda) $  and $s_{jl}(\lambda) $ are some numerical coefficients.

  These formulas resemble of course asymptotic formulas for eigenfunctions of differential operators
      \begin{equation}
B=    D^n+ \sum_{m=0}^{n-1} b_m (x) D^m, \q  \q D   = D_{x}   =- id/d x,
 \label{eq:VB}\end{equation}
in the space $L^{2} (\Bbb R)$ with  the  coefficients $b_m (x)$, $m=0,1,\ldots, n-1$,   decaying sufficiently rapidly at infinity. The terms $t^{-1/2}   e^{  \pm i \omega (t,k)}$ play the role of the functions $  e^{  \pm i kx}$ for differential operators. They are also similar to  asymptotic formulas of \cite{Y2} for eigenfunctions of the perturbed Carleman  operator; in this case $\omega (t,k)=k \ln t$ and $\lambda(k)= \pi (\cosh (\pi k))^{-1}$. 

Relations \e{eq:Sodd2H}, \e{eq:So-}  show that, for $n$ odd, a wave coming from zero (from infinity) cannot penetrate to infinity (zero) so that there is the complete reflection in this case. For $n$ even, $s_{11}(\lambda) $ and $s_{22}(\lambda) $ 
in relations \e{eq:S1H}, \e{eq:S2H} are naturally interpreted as the reflection coefficients while $s_{21}(\lambda) $ and $s_{12}(\lambda) $ are   interpreted as the transmission coefficients. We show  that  
$s_{21}(\lambda) = s_{12}(\lambda) $.

As could be expected the  coefficients $s(\lambda) $ in \e{eq:Sodd2H} and $s_{jl}(\lambda) $ in \e{eq:S1H}, \e{eq:S2H} are elements of appropriate scattering matrices $S(\lambda)$. They are unitary; in particular,  $|s (\lambda)|=1 $ in \e{eq:Sodd2H}. However $S(\lambda)$ are the scattering matrices for some pair of auxiliary differential operators and not for Hankel operators themselves. In fact, there is no ``unperturbed''  Hankel operator $H_{0}$ such that $S(\lambda)$ is the scattering matrix for the pair $H_{0}, H$.

\medskip
    
     {\bf 1.3.}
     Our approach relies on a reduction of Hankel operators with kernels  \e{eq:LOG},  \e{eq:LOG4}  to differential operators.
   It was shown in \cite{Yf1} that such Hankel operators  are unitarily equivalent
 (essentially, by the Mellin transform) to the  operators
  \begin{equation}
A= v Q (D_{\xi} ) v ,\q D_{\xi}   = -i d /d \xi,
\label{eq:LOGz}\end{equation}
 in the space $L^2 ({\Bbb R}) $.
Here $v$ is the operator of multiplication by the universal function 
  \begin{equation}
v (\xi )=\frac{\sqrt{\pi}} {\sqrt{\cosh (\pi \xi)}}
 \label{eq:LOGZ1}\end{equation}
and the real polynomial  
 \begin{equation}
 Q(X)=  \sum_{ m=0}^n q_m X^m 
\label{eq:LO4}\end{equation}
 is determined by $P (X )$. The polynomials $P (X)$ and $Q (X )$ have the same degree,   and their coefficients   are linked by an explicit formula (see formula \e{eq:LOG5a} below); in particular, $ q_{n}   =1$ if $  p_n =1$. If $n=0$, then $Q(X )=P( X)=1$ so that $A=v^2$. This yields the familiar diagonalization of the Carleman operator.
 
 Actually, it was shown in \cite{Y,Yafaev3} that every Hankel operator is  unitarily equivalent to the pseudodifferential operator  \e{eq:LOGz}  with the function   $Q(X)$ determined by the kernel $h(t)$ and called the sign-function of $H$
 in \cite{Y}.
 This terminology is explained by the fact that a Hankel operator $H\geq 0$ if and only if the function   $Q(X)\geq 0$.  In general, 
 $Q(X)$ is a distribution. So this paper is devoted to a study of the case when the sign-function is an arbitrary polynomial.

    Observe that the highest order term of the operator $A$ equals $a_n (\xi) D_{\xi}^n $
 where $a_n (\xi)= v^2 (\xi) $ tends to zero (exponentially) as $|\xi|\to\infty$.
                    Apparently such differential operators  were never studied before, and we are led to fill in this gap.
Studying   differential operators \e{eq:LOGz}  we do not make specific assumption \e{eq:LOGZ1} and consider sufficiently arbitrary real  functions $v(\xi)$  tending to zero as $|\xi | \to\infty$.
     The essential spectrum of differential operators \e{eq:LOGz} 
 was  described in \cite{Yf1}    where it was shown that $\spec  (A)={\Bbb R}$ if $n$ is odd, and $\spec_{\rm ess} (A)= [0, \infty)$ if $n$ is even. 
   The last result should be compared with the fact that $\spec_{\rm ess} (A)= [\min Q( X), \infty)$  if $v(\xi)=1$. Thus, even in this relatively simple question,  the degeneracy of $v(\xi )$ at infinity significantly changes spectral properties of  differential operators $A$.
        Here we study the detailed spectral structure, in particular, the absolutely continuous spectrum,  of differential operators \e{eq:LOGz} and hence of Hankel operators with kernels \e{eq:LOG}. 
  
We show that differential operators \e{eq:LOGz}
    can be reduced by an explicit unitary transformation $\mathrm {L}$ (the generalized  Liouville transformation)  to standard differential operators. Set
    \[
 (\mathrm {L} u) (\xi)=  x'(\xi)^{1/2} u( x(\xi))
 \]
 where the variables $x$ and $\xi$ are linked by  the relation
\[
x= x(\xi)= \int_0^\xi v(\eta)^{-2/n} d\eta
\]
 so that $x'(\xi)=v(\xi)^{-2/n}$. Then the operator
$
B=\mathrm {L}^* A \mathrm {L}
$
 is given by the formula     \e{eq:VB}.
 Our crucial observation is that   the  coefficients $b_m (x)$, $m=0,1,\ldots, n-1$, of this differential operator   decay at infinity.    
 
     It is noteworthy that if $v(\xi)$  (for example,  function  \e{eq:LOGZ1}) tends to zero exponentially as $|\xi|\to\infty$, then the coefficients (except $b_{n-1} (x)$ which can be removed by a gauge transformation) of the operator $B$ decay faster than $|x|^{-1}$ as $|x|\to\infty$. On the contrary, for slower decay of $v(\xi)$, the   coefficients of the operator $B$   decay slower or as   $|x|^{-1}$. We say that in these cases    the operator $B$ has short-range or long-range coefficients, respectively.
   Thus, somewhat counter-intuitively,  a stronger     degeneracy of the operator $A$ yields  better properties of the operator $B$.
   
   \medskip

 {\bf 1.4.}
 Let us briefly  describe  the structure of the paper.
 We first consider differential operators \e{eq:VB} with decaying coefficients $b_m (x)$. By a gauge transformation, we can obtain the operator $B$ of the  same structure but with $b_{n-1}(x)=0$. Thus for $n=1$, the operator $B$ reduces to the operator $D$ so that the problem is trivial. For $n=2$, the spectral analysis of differential operators \e{eq:VB} is a  very well developed machinery both for short- and long-range coefficients  $b_m (x)$. As far as the case $n>2$ is concerned, we note that the methods of functional analysis work equally well for all $n$. On the contrary, there is a substantial difference in application of specific methods of differential equations. This difference is particularly important in the long-range case and, for $n>2$, we have to carry out an analysis which was not available in the literature. 
 The cases of short- and long-range coefficients are considered in Sections~2 and 3, respectively.
 
 In Section~4 we show that degenerate differential operators  \e{eq:LOGz}  can be reduced by   the generalized  Liouville transformation  to  operators \e{eq:VB} with the coefficients $b_{m}(x)$ decaying at infinity. This leads to new spectral results for operators \e{eq:LOGz}. We emphasize that for function  \e{eq:LOGZ1}, this reduction together with the gauge transformation yields the operator $B$ with short-range coefficients.

We return to Hankel operators   $H$     in Section~5. Here we proceed from the fact that Hankel operators with kernels  \e{eq:LOG}  are unitarily equivalent to degenerate differential operators $A$ given by \e{eq:LOGz}. Therefore the results of  Section~4 allow us to prove Theorem~\ref{SpThHa}.
Spectral results (i), (iii) and (iv) of Theorem~\ref{SpThHa} are direct consequences of the corresponding results for the operators $A$ and hence of the same results for the differential  operators $B$ with short-range coefficients. On the  contrary, the limiting absorption principle (statement (ii)) for the operator $H$ does not follow from the corresponding statement
for $A$. It requires estimates on eigenfunctions $\theta(t,k)$ of  the operator $H$   which are obtained as consequences of asymptotic formulas \e{eq:Sodd2H}, \e{eq:So-}  or \e{eq:S1H}, \e{eq:S2H}. These formulas are proven     in Section~6 where the evolution operator $e^{-iHT}$ is also studied for $T\to\pm\infty$. This is technically the most difficult part of the 
paper.  Finally,  in Section~7, we study the asymptotic behavior of the unitary group $e^{-iHT}$ as $T\to\pm\infty$.
 
 %reformulate a part of our results on the integral operators \e{eq:H1}  in the framework of the discrete realization of Hankel operators.

\medskip
 
 xxxxxxxxx   
    
 {\bf 1.5.}
   Let us introduce some standard
 notation. 
 %We first recall that ${\Bbb T}$ is the unit circle in the complex plane and ${\Bbb Z}_{+} $ is the set of all nonnegative integers.
 We denote by $\Phi$,
 \[
(\Phi f) (k)=  (2\pi)^{-1/2} \int_{-\infty}^\infty f(x) e^{ -i x k} dx,
\]
 the Fourier transform.        We often use the same notation for a function and the operator of multiplication by this function.       
     We denote by   $ {\sf H}^k ({\cal J})$  the Sobolev space of functions defined on an interval ${\cal J}\subset {\Bbb R}$; 
     $ C_{0}^k({\cal J})$ is the class of $k$-times continuously differentiable functions with compact supports in ${\cal J}$.
   The letters $c$ and $C$ (sometimes with indices) denote various positive constants whose precise values are inessential.

     %in particular,         ${\sf H}^K={\sf H}^K({\Bbb R})$.

 %  The Dirac function is standardly denoted $\d(\cdot)$; $\d_{n,m}$ is the Kronecker symbol, i.e., $\d_{n,n}=1$  and $\d_{n,m}= 0$  if $n\neq m$.   

   %************************************************************
\section{Differential operators of arbitrary order}  
%***********************************************************

 {\bf 2.1.}
  Here we consider  differential operators $B$ defined by equality \e{eq:VB}.
 We always suppose that    operators $B$ are symmetric on $C_{0}^\infty ({\Bbb R})$, that is, 
   \[
 \sum_{m=0}^{n-1} b_m (x) D^m =  \sum_{m=0}^{n-1}    D^m \ov{b_m (x) }
\]
% \label{eq:SAd}\end{equation}
 where the derivatives of $b_m (x) $ are understood   in the sense of distributions.
  This leads to certain algebraic relations for
   the functions $b_m (x) $ which we assume to be satisfied; in particular, $b_{n-1} (x) =\ov{b_{n-1} (x) }$.  If   the coefficients $b_m (x) $, $m=0,1, \ldots, n-1$,  are bounded, then the   operator $B$
 is self-adjoint on    
    the Sobolev class ${\sf H}^n ({\Bbb R})$. If, moreover, $b_m (x) \to 0$ as $|x|\to\infty$, then  the perturbation $V=B-B_{0}$ of the ``free" operator $B_{0}=D^n$ is compact relative to  $B_{0}$. It follows that the essential spectra of the  operators $B_{0} $ and $B$ are the same. Thus the essential spectrum of the operator $B$ coincides with $\Bbb R$ if $n$ is odd and it coincides with $[0,\infty)$ if $n\geq 2$ is even.

Suppose now that the coefficients of the operator  \e{eq:VB} are short-range.
The following   assertion contains basic results of  spectral analysis  of such differential operators.

  \begin{theorem}\label{SpTh}
  Assume that
    \begin{equation}
| b_m  (x)| \leq C (1+|x|)^{-\rho }, \q \rho >1,\q m=0,1, \ldots, n-1,
 \label{eq:SR}\end{equation}
and let $B $ be the self-adjoint   operator defined by differential  expression \e{eq:VB} on the Sobolev class ${\sf H}^n ({\Bbb R})$.  Then:
\begin{enumerate}[\rm(i)] 
  \item
   The spectrum of the operator $B$ is absolutely continuous except   eigenvalues that  may accumulate to zero and  infinity only.  
 \item
  The limiting absorption principle holds, that is, for any $\sigma>1/2$, the operator-valued  function 
  \[
  \la x\ra^{-\sigma} ( B-z)^{-1} \la x\ra^{-\sigma}, \q \Im z \neq 0 ,
  \]
 is H\"older continuous with exponent $\alpha<\sigma-1/2$ $($and $\alpha<1)$  up to the real axis,  except the eigenvalues of the operator $B$ and the point zero.
\item
     The 
 absolutely continuous  spectrum of the operator $B$ covers $\Bbb R$ and is simple for $n$ odd. It coincides with $[0,\infty)$ and has multiplicity $2$  for $n$ even.  
 \item
If $n$ is odd, then the multiplicities of eigenvalues are bounded by $(n-1)/2$. If $n$ is even, then the multiplicities of positive eigenvalues are bounded by $n/2-1$, and the multiplicities of negative eigenvalues are bounded by $n/2$.
% The point zero is not an eigenvalue.
   \end{enumerate} \end{theorem}
    
    Theorem~\ref{SpTh} is a standard result of scattering theory for the operators 
     $B_{0}$ and $B$. Below we describe briefly main steps of its proof relying on the smooth approach in abstract scattering theory. Note that for 
     $n\geq 3$  specific methods of ordinary differential equations are not convenient so that the case $n\geq 3$ is closer to multidimensional problems than to   operators  \e{eq:VB} for $n=2$.
 A proof of parts (i) and (ii) can also be obtained by the Mourre method which is discussed in Section~3. As far as the inverse scattering problem for  differential operators of arbitrary order is concerned, we refer to the book \cite{BDT}.

The smooth approach requires a preliminary study
of the operator $B_{0}=  D^n$. Its spectrum   is absolutely continuous. For odd $n$, it is simple and   coincides with $\Bbb R$. For even $n$, it  has multiplicity two  and coincides with $[0,\infty)$.   Using the Fourier transform it  is easy to calculate the integral kernel $R_{0}(x,y;z)$ of the resolvent  $R_{0}(z)= (B_{0}-z)^{-1}$, $\Im z\neq 0$, of the operator $B_0 $. Let $\zeta_{j}$, $ j=1, \ldots, n$,  be the solutions of the equation $  \zeta^n=z$. Then
\begin{equation}
\begin{split}
 R_0(x,y; z) &=   i n ^{-1}\sum_{\Im \zeta_{j}>0} \zeta_{j}^{-n+1} e^{i \zeta_{j} (x-y)}, \q x\geq y,
 \\
  R_0(x,y; z)& = -  i n  ^{-1}\sum_{\Im \zeta_{j}<0} \zeta_{j}^{-n+1} e^{i \zeta_{j} (x-y)}, \q x\leq y .
\end{split}
\label{eq:R}\end{equation}
Since the operator $B_{0}$ commutes (anticommutes) with the complex conjugation for even (odd) $n$, we have the identities
%\begin{equation}
\[
\begin{split}
 R_0(x,y; \bar{z}) &=   \ov{R_0(x,y; z)}, \q {\rm even}\q n,
 \\
  R_0(x,y; \bar{z})& = -  \ov{R_0(x,y; -z)} ,  \q {\rm odd}\q n,
\end{split}
\]
%\label{eq:Rc}\end{equation}
which follow also from explicit formulas \e{eq:R}.  
Obviously, for fixed $x,y$, the analytic function $ R_0(x,y; z) $ of $z$, $\Im z\neq 0$, is continuous up to the real axis with   exception of the point $z=0$. Moreover, for even $n$,  this function  is actually analytic in the complex plane cut along $[0,\infty)$. Formulas \e{eq:R} imply also that the operator-valued functions 
\begin{equation}
  \la x\ra^{-\sigma} D^p R_{0}(z) \la x\ra^{-\sigma}, \q \sigma>1/2, \q p=0,\ldots, n-1,
\label{eq:Rt}\end{equation}
  of $z$ possess the same properties of analyticity and continuity in the Hilbert-Schmidt norm. 
  
%  We note also some asymptotic formulas. Suppose that $k\in{\Bbb R}\setminus\{ 0\}$ and $f\in  L^2 ({\Bbb R})\cap L^1 ({\Bbb R})$. If $n$ is odd, it follows from \e{eq:R} that
%  \begin{equation} ( R_{0}(k^n + i0) f)(x)=  i n ^{-1} k^{-n+1} e^{i  k x} \int_{-\infty}^\infty e^{-i  k y} f(y) dy, \q x \to +\infty,
%\label{eq:Rodd}\end{equation} and $( R_{0}(k^n + i0) f)(x)\to 0$ as $x\to-\infty$. Similarly, if $n$ is even, then
%  \begin{equation} ( R_{0}(k^n + i0) f)(x)=  i n ^{-1} |k|^{-n+1} e^{\pm i | k |x} \int_{-\infty}^\infty e^{\mp i  k y}  f(y)  dy, \q x \to \pm\infty.\label{eq:Reven}\end{equation}
  
  To obtain similar information on the resolvent $R (z) =(B -z)^{-1}$  of the operator $B$, we regard
   the   resolvent identity
\begin{equation}
R(z)= R_{0}(z)-R  (z)VR_{0} (z), \q \Im z\neq 0, \q V=B -B_{0},
\label{eq:res}\end{equation}
 as the Fredholm equation for $R (z)$. Put $G=  \la x\ra^{-\rho/2} $, $G_{0}= G^{-1} V$. Then it follows from \e{eq:res} that
 \[
G D^{p} R(z) G = G D^{p} R_{0}(z)G \big(I + G_{0} R_{0}(z)G \big)^{-1} ,\q p=0,1,\ldots, n-1,
\]
%\label{eq:res1}\end{equation}
where the inverse operator in the right-hand side exists in view of the self-adjointness of the operator $B$.  Let $\cal N \subset {\Bbb R}\setminus\{0\}$ for odd $n$ and $\cal N\subset {\Bbb R}_{+}\setminus\{0\}$ for even $n$ be the set of $\lambda$  where at least one of two homogeneous equations
\[
f + G_{0} R_{0}(\lambda\pm i 0)Gf=0
\]
% \label{eq:res2}\end{equation}
has a nontrivial solution $f\in L^2 ({\Bbb R})$.  
 According to the analytic Fredholm alternative (see, e.g., \cite{YaMSC}, Theorem~1.8.2) the set $\cal N$ is closed and has the Lebesgue measure zero. On the complement of the set $\cal N$, the spectrum of the operator $B$ is absolutely continuous and the operator-valued function $G D^{p} R(z) G $ is continuous (in the Hilbert-Schmidt norm) as $z$ approaches the cut along the continuous spectrum. 
 
 The next step is to   prove  that $\cal N$ consists of eigenvalues of the operator $B$, and these eigenvalues may accumulate to the point zero only. In particular,  the operator $B$ does not have the  singular continuous spectrum.    For the proof, one can use
  the scheme of S.~Agmon \cite{AgC} which  significantly
    simplifies (see \cite{Y4ord})  in the one-dimensional case. Thus we obtain statements (i) and (ii) of Theorem~\ref{SpTh}. In fact, instead of (ii) we   prove a slightly stronger assertion.

  \begin{proposition}\label{SpThDD}
  Under assumption
 \e{eq:SR}   the operator-valued  function 
  \[
  \la x\ra^{-\sigma} D^{p}( B-z)^{-1} \la x\ra^{-\sigma}, \q \Im z \neq 0 , \q p=0,1, \ldots, n-1, \q \sigma>1/2,
  \]
 is H\"older continuous with exponent $\alpha<\sigma-1/2$ $($and $\alpha<1)$  up to the real axis,  except the eigenvalues of the operator $B$ and the point zero.
   \end{proposition}

% Statement (ii) relies on the results of scattering theory.

% For a self-adjoint operator  $B$ in a Hilbert space $\cal H$, 
% let  ${\cal H}^{(p)}$ be the subspace of $\cal H$ spanned by all eigenvectors of the operator $B$;
% ${\cal H}^{(c)}$ is its orthogonal complement and ${\cal H}^{(c)}$ is the absolutely continuous subspace of $B$. We denote  by $P^{(ac)}$ be the orthogonal on ${\cal H}^{(c)}$. 
%$P^{(c)}$ is the orthogonal projection onto ${\cal H}^{(c)}$;  $B^{(c)}$ is the restriction of $B$ onto ${\cal H}^{(c)}$. It follows from (i) that the operator $B^{(c)}$ is absolutely continuous.

   \medskip
   
   %  \begin{equation}
% W_\pm=W_\pm(B , B_{0}  )  =\slim_{t\to\pm\infty}e^{i B t}  e^{-i B_{0} t}. 
% \label{eq:WO}\end{equation}
 
 {\bf 2.2.}
 Let us now discuss statement (iii).
First, we briefly recall the definition of wave operators for a pair of self-adjoint operators  $B_{0}$, $B$ acting in a Hilbert space  $\cal H$. Let $P^{(ac)}$ be the orthogonal projection on the absolutely continuous subspace ${\cal H}^{(ac)}$ of the operator $B$, and let $B^{(ac)}$ be the restriction of $B$ on ${\cal H}^{(ac)}$. Similar objects for the operator $B_{0}$ will be endowed with the index $``0"$.
  The wave operators for a pair  $B_{0}$, $B$ and a bounded operator  $J  $ (``identification")  are   defined as strong limits
     \begin{equation}
W_{\pm}= W_{\pm} (B,B_{0};J)= \slim_{T\to\pm\infty}e^{iB T} J e^{- iB_{0}T}P_{0}^{(ac)}. 
\label{eq:WLR}\end{equation}
 Under the assumption of their existence,  the wave operators   \e{eq:WLR} enjoy the intertwining property 
 $B W_{\pm} =W_{\pm} B_{0}$.
We also consider the   wave operators 
    \begin{equation}
  W_{\pm} (B_{0}, B;J^*)= \slim_{T\to\pm\infty}e^{iB_{0}T} J^* e^{- iB T} P^{(ac)} 
\label{eq:WLRad}\end{equation}
 for the pair $B $, $B_{0}$ and the  ``identification" $J^*$.
  If both limits  \e{eq:WLR}  and \e{eq:WLRad}   exist, then they
   are adjoint to each other.
   
 Let $I$  be the identity operator.  In the important particular case $J=I$  the operator $W_{\pm} (B,B_{0}):=W_{\pm} (B,B_{0}; I)$ is isometric on ${\cal H}_{0}^{(ac)}$ and the existence of the wave operator $W_{\pm} (B_{0}, B)$ is equivalent to the relation
   \[
   \Ran (W_{\pm} (B,B_{0}))={\cal H}^{(ac)}
   \] 
   known as the completeness of $W_{\pm} (B, B_{0})$. This relation implies that the operators $B_{0}^{(ac)}$ and $ B^{(ac)}$ are unitarily equivalent.
 
For the proof of the existence of   limits \e{eq:WLR} and \e{eq:WLRad}, we rely on the T.~Kato theory of smooth perturbations
(see, e.g., \S XIII.7.C of the book \cite{RS} or  \S 4.5 of the book \cite{YaMSC}). We use the following

\begin{definition}\label{Kato}
 A $B$-bounded operator $G$ is called $B$-smooth  (in the sense of Kato)  if for some open set $\Omega\subset \spec (B)$ of full measure  (that is,  the Lebesgue measure $|\spec (B)\setminus \Omega|=0$) and every compact interval $X\subset \Omega$, we have
   \begin{equation}
\sup_{\lambda\in X, \varepsilon \neq 0} \| G \big( R(\lambda+i\varepsilon) -  R(\lambda-i\varepsilon)\big) G^*\| <\infty. 
\label{eq:smooth}\end{equation}
 \end{definition}

In particular,   condition \e{eq:smooth} is satisfied if the operator-valued function $G R(z) G^*$, or more generally 
$(G R(i)) R(z) (G R(i))^*$, is continuous as $z$ approaches the set $\Omega$ (from the upper and lower half-planes).
We need the following result.

\begin{proposition}\label{smooth}
Suppose that
 \[
 BJ-JB_{0}= G^* G_{0}
 \]
%\label{eq:smooth1}\end{equation}
where the operators  $G_{0}$ and $G$ are $B_{0}$- and $B$-smooth, respectively. Then the wave operators \e{eq:WLR}  and \e{eq:WLRad} exist.
 \end{proposition}

%  In application to the operators $B_{0}$, $B$ defined by \e{eq:VB} we have the following result (see, e.g., Theorem~4.5.6 in \cite{YaMSC}).  Now a priori we only assume that the coefficients  $b_{m}(x)$  are bounded and the symmetricity condition \e{eqSad} is satisfied.

Let us come back to differential operators $B_{0}= D^n$ and $B$ defined by formula \e{eq:VB}. Now $J=I$,
$P_{0}^{(ac)}=I$,   $G=G^*=\la x\ra^{-\rho/2}$ and $G_{0}= G^{-1}V$. 
It follows from the continuity of the operator-valued function \e{eq:Rt} that the operator $G_{0}$ is $B_{0}$-smooth; the corresponding open set $\Omega_{0}$ of full measure is 
$ {\Bbb R}_{+} $ for $n$ even and  $\Omega_{0} = {\Bbb R}\setminus \{0\} $ for $n$ odd.
According to
part (ii) of Theorem~\ref{SpTh} the same result is true for the pair of the operators $G$ and   $B$; one only has to remove its point spectrum
$\spec_{p}(B)$. Thus the operator $G$ is $B$-smooth; the corresponding set $\Omega$ is given by the relations
 \begin{equation}
\Omega={\Bbb R}_{+}\setminus\spec_{p}(B)\; {\rm    for} \; n \; {\rm  even }
\; {\rm    and}\;  \Omega={\Bbb R}\setminus (\{0\}\cup \spec_{p}(B)) \;{\rm  for }\;  n \;  {\rm odd}.
\label{eq:omega}\end{equation}
Therefore Proposition~\ref{smooth}  yields the following result.

% Let the operator $J$ be defined by formulas    \e{eq:Ga} and  \e{eq:Ga1}.   

\begin{theorem}\label{SpThSRa}
Under the assumptions of Theorem~\ref{SpTh}  the wave operators $ W_{\pm} = W_{\pm} (B,B_{0})$  exist.     
  They   
 are isometric and complete, that is, their ranges $\ran ( W_\pm )  ={\cal H}^{(ac)}$. The intertwining property 
$ B W_\pm=W_\pm B_{0} $
holds. The operators $B^{(ac)}$ and $B_{0}$ are unitarily equivalent.
 \end{theorem}
  
  Of course Theorem~\ref{SpThSRa} implies that the absolutely continuous spectrum of the operator $B$ is the same as that of $B_{0}$. This concludes the proof of part (iii) of Theorem~\ref{SpTh}.

  \begin{remark}\label{SCTh}
 The existence of the wave operators $W_{\pm}=W_{\pm} (B,B_{0})$ is   a simple result that does not require the theory of smooth perturbations. It can be obtained by an elementary Cook's method (see, e.g.,  \S 1.4 of the book \cite{Ya}).  The existence of   $W_{\pm} $ entails that the restriction of $B$
 on the subspace $\Ran(W_{\pm})$ is unitarily equivalent to the operator $B_{0}$. This fact can be combined
 with  the Weyl-Titchmarsh-Kodaira theory (see, e.g.,  the books \cite{CoLe, Nai}) which implies that the multiplicity of the spectrum of the operator $B$ does not exceed $n$. For $n=2$ (but not for larger $n$) this yields   the completeness of the   wave operators $W_{\pm}$. 
  \end{remark}
    
  \begin{remark}\label{SCThX}
   The existence and   completeness of   $W_{\pm}$  are also  consequences of the trace class scattering theory
  (see, e.g.,  Chapter 6 of \cite{YaMSC}).
  \end{remark}
  
  Finally, we discuss part (iv) of Theorem~\ref{SpTh}. It can be easily checked (see, e.g.,  \cite{CoLe, Y4ord}) that for every $\lambda\in {\Bbb R}\setminus\{ 0\}$ the differential equation $Bf =\lambda f$ has solutions $f_{j}$, $j=1,\ldots, n$, such that
    \begin{equation}
f_{j} (x,\lambda)=e^{i k_{j}  x} (1+ o(1)), \q k_{j}^{n}=\lambda, 
\label{eq:EigF}\end{equation}
as $x\to +\infty$. These solutions are linearly independent and     $f_{j} \in L^{2} ({\Bbb R}_{+})$ if and only if $\Im \kappa_{j} > 0$. For $n$ odd, the number of such $\kappa_{j}$
  equals $(n-1)/2$. For $n$ even, $\Im \kappa_{j} > 0$ for $n/2 -1$ values of $j$ if $\lambda>0$ and 
 for $n/2 $ values of $j$ if $\lambda<0$. So we obtain the upper bound on the multiplicities of eigenvalues stated in  part (iv). This concludes the proof of Theorem~\ref{SpTh}.

   Note that we have distinguished solutions $f_{j} (x,\lambda)$ by their asymptotics \e{eq:EigF}  for $x\to\infty$ only for definiteness. The same arguments work for $x\to-\infty$.    
 
 \medskip

   {\bf 2.3.}
   The construction of the wave operators is intimately related with eigenfunction expansions.   This is the classical stuff for $n=2$; we refer to the paper \cite{F} by L.~D.~Faddeed or  the book \cite{Ya}, Chapters~4 and 5. 
 For all $n\geq 3$, the construction below is probably not explicitly written in the literature, but it is essentially the same as for the particular case $n=4$ discussed in \cite{Y4ord}. We note that in contrast to the case $n=2$ when one can use Volterra integral equations, for $n\geq 3$ one is obliged to work with  Fredholm equations. From this point of view, the case $n\geq 3$ is closer to multidimensional problems than to  operators  \e{eq:VB} for $n=2$.
 
 The ``free" operator $B_{0}$ can of course be diagonalized by the Fourier transform  $\Psi_{0}$, 
  \[
({\Psi }_{0} f)(k)= (2\pi)^{-1/2}  
 \int_{-\infty}^\infty e^{-ikx}f(x) dx ,
 \]
 %\label{eq:Four}\end{equation}
  that is,  $\Psi_{0}B_{0}=\Lambda\Psi_{0}$ where $\Lambda$ is the operator of multiplication by the function $k^n$ in the space $L^2 ({\Bbb R})$. It means that 
 $\psi_{0} (x,k)=e^{ikx}$ where $k\in {\Bbb R}$ is a complete set of eigenfunctions of the operator $B_{0}$.   Let us define  eigenfunctions of the continuous spectrum of the operator $B$
  by the relation
  \begin{equation}
\psi_{\pm} (k)= \psi_{0}(k) -R (k^n \mp i 0) V \psi_{0} (k) , \q \lambda= k^n\in \Omega.
\label{eq:WF}\end{equation} 
By virtue of the limiting absorption principle (Proposition~\ref{SpThDD})  the right-hand side here is correctly defined and $\la x \ra^{-\sigma} D^{p}\psi_{\pm}\in L^2 (\Bbb R)$ for any $\sigma>1/2$ and $p=0,1,\ldots, n-1$.  
The resolvent identity \e{eq:res} implies that the functions
$\psi _{ \pm} (x, k)$ satisfy    the Lippmann-Schwinger equation
\begin{equation}
\psi_{ \pm}(k)= \psi_{0} (k) -R_{0}(k^n \mp  i 0) V \psi_{ \pm}(k)  .
\label{eq:WF2}\end{equation}
Since, by \e{eq:R}, the integral kernel of the operator $D^{p}R_{0}(k^n \mp  i 0) $ is a bounded function, we have the estimate
  \begin{multline*}
| (D^{p} R_{0}(k^n \mp  i 0) V \psi_{ \pm}(k) )(x)|\leq C\int_{-\infty}^{\infty} | (V \psi_{ \pm})(x,k) | dx
\\
\leq C\sum_{m=0}^{n-1}
\| b_{m }\la x \ra^{\rho/2}\| \|  \la x \ra^{-\rho/2} D^{m}\psi_{ \pm}(k)\| .
  \end{multline*}
  In view of \e{eq:WF2}, this yields the estimate
  \begin{equation}
| D^{p}\psi_{ \pm}(x,k)| \leq C (k) , \q \lambda= k^n\in \Omega, \q p=0,1, \ldots, n-1.
\label{eq:WFDD}\end{equation}
   According to \e{eq:WF} the functions $\psi_{\pm} (k)$ satisfy also  the differential equation $B \psi_{\pm} (k)=k^n \psi_{\pm} (k)$, that is,
   \begin{equation}
i^{-n}\psi_{\pm}^{(n)}(x, k)+ \sum_{m=0}^{n-1} i^{-m} b_m (x) \psi_{\pm}^{(m)}(x , k)= k^n\psi_{\pm}(x, k).
 \label{eq:DE}\end{equation}
Of course for even $n$, the functions $\psi_{\pm}(x, k)$ and $\psi_{\pm}(x, -k)$  satisfy the same equation \e{eq:DE} while these equations are different if $n$ is odd. 
 
 For  the operator  $B$, there  exist two natural  diagonalizing transformations  denoted   $ \Psi_\pm$.
They are constructed in terms of the functions $\psi_{\pm}(x, k)$   by the formula
 \begin{equation}
({\Psi }_{\pm} f)(k)= (2\pi)^{-1/2}  
 \int_{-\infty}^\infty \overline{\psi_{\pm}(x, k)}f(x) dx , \q f\in C_{0}^\infty ({\Bbb R} ), \q  \lambda = k^n \in \Omega,
 \label{eq:WFsr1}\end{equation}
(actually it suffices to assume that $\la x \ra^{\sigma} f\in L^2 (\Bbb R)$ for some  $\sigma>1/2$).  
 These mappings   extend by continuity to 
bounded operators on $ L^{2}({\Bbb R} )$; they satisfy the relation 
\begin{equation}
\Psi_{\pm} \Psi_{\pm}^* =  I, \q \Psi_{\pm}^*\Psi_{\pm}=P^{(ac)} 
 \label{eq:F1}\end{equation}
 and diagonalize $B$, that is,
 \begin{equation}
\Psi_{\pm} B =\Lambda \Psi_\pm.
 \label{eq:F2}\end{equation} 
 The wave operators $W_{\pm}  = W_{\pm} (B,B_{0})$ and the operators \e{eq:WFsr1} are linked by the formula
  \begin{equation}
W_{\pm} =  \Psi_{\pm}^*  \Psi_{0}.
 \label{eq:WTD}\end{equation} 
 It follows from formulas \e{eq:F1} and \e{eq:F2} that the spectral family $E(\lambda)$ of the operator $B$ satisfies the identity 
 \begin{equation} 
\frac{d(E (\lambda)f,f)}{d\lambda}=\sum_{k^n=\lambda}\big|\int_{-\infty}^\infty   \ov{{ \psi}_{\pm} (x,k)} f( x) d x\big|^2, \q f\in C_{0}^\infty ({\Bbb R} ), \q  \lambda = k^n \in \Omega,
\label{eq:est2B}\end{equation}
holds.
Here the sum consists of the one term if $n$ is odd; if $n$ is even it consists of two terms for $\lambda> 0$, and it is empty for $\lambda<0$.

   The scattering operator 
 \begin{equation}
{\cal S}:=W_{+}^* W_{-}
 \label{eq:ScOp}\end{equation}
  commutes  with the operator $B_{0}$ and  according to \e{eq:WTD}
\begin{equation}
\Psi_{0}{\cal S} \Psi_{0}^*= \Psi_{+} \Psi_{-}^*.
\label{eq:S}\end{equation}

Let us summarize the results described  above.

 \begin{theorem}\label{SpThSR}
 Let estimates \e{eq:SR} be satisfied.
  Define the functions $\psi_{\pm} (x,k)$  by formula  \e{eq:WF} and the mappings $\Psi_{\pm}  $    by formula  \e{eq:WFsr1}. Then the relations \e{eq:F1}, \e{eq:F2} and  \e{eq:WTD}  hold. The scattering operator \e{eq:ScOp} satisfies identity \e{eq:S}. It is a unitary operator.
    \end{theorem}

Recall that the operator 
 \e{eq:ScOp} 
  commutes  with the operator $B_{0}$. Therefore,
for odd $n$, the operator $\wh{\cal S} = \Psi_{0}{\cal S}\Psi_{0}^{*}$   acts in $L^{2 }({\Bbb R})$ as the multiplication by the function $s(\lambda)$, that is,
\begin{equation}
( \wh{\cal S} g)(k)=s(\lambda) g(k), \q \lambda=k^n\in \Omega \subset{\Bbb R}.
\label{eq:S1odd}\end{equation}
%  where the function $s(\lambda)$   is the same as in  \e{eq:WFSodd2}. 
For  even $n$, we set $(\mathrm {Y} g)(k)=(g(k),g(-k))^{\top}$, $k>0$; obviously, the operator $\mathrm {Y}: L^{2}({\Bbb R})\to L^{2}({\Bbb R}_{+})\otimes {\Bbb C}^{2}$ is unitary. The operator $\mathrm {Y}\wh{\cal S} \mathrm {Y}^{*}$   acts in $L^{2}({\Bbb R}_{+})\otimes {\Bbb C}^{2}$ as the multiplication by the the matrix-valued function $S(\lambda)$, that is,
%\begin{equation}
\[
   \left(\begin{array}{cc}
(\wh {\cal S} g)(k)
\\
(\wh {\cal S} g)(-k)
\end{array} \right)  = S(\lambda) 
  \left(\begin{array}{cc}
g (k)
\\
g(-k)
\end{array} \right), \q  \lambda=k^n\in \Omega\subset{\Bbb R}_{+}, \q k>0, 
\]
%\label{eq:S1}\end{equation}
where   
  \begin{equation}
S(\lambda) =\left(\begin{array}{cc}
s_{11}(\lambda)&s_{12}(\lambda)
\\
s_{21}(\lambda)&s_{22}(\lambda)  \end{array} \right) 
\label{eq:SM}\end{equation} 
 is the $2\times 2$ matrix.
The scalar function $s(\lambda)$ for $n$ odd and the matrix-valued function $S(\lambda)$ 
for $n$ even  are 
known as the scattering matrices.  The  scattering matrices are well defined for $\lambda\in\Omega$. They are unitary and depend continuously on $\lambda\in\Omega$. 
   
 It follows from formula \e{eq:S} that $\Psi_{+}^*=\Psi_{-}^* \wh{\cal S}^*$. Therefore using definition
 \e{eq:WFsr1} of the operators $\Psi_\pm$,  we see that $\psi_{+}(x,k) = \ov{s(\lambda)}\psi_{-}(x,k) $ if $n$ is odd and
 \begin{equation}
 \left(\begin{array}{cc}
\psi_{+}(x,k)
\\
\psi_{+}(x,-k)
\end{array} \right)  =\left(\begin{array}{cc}
\ov{s_{11}(\lambda)}&\ov{s_{12}(\lambda)}
\\
\ov{s_{21}(\lambda)}& \ov{s_{22}(\lambda)}  \end{array} \right)
  \left(\begin{array}{cc}
\psi_{-}(x,k)
\\
\psi_{-}(x,-k)
\end{array} \right)
\label{eq:SMev}\end{equation} 
if $n$ is even and $k>0$. 
     
    \medskip
 
  {\bf 2.4.} 
   Let us find the asymptotic behavior of the functions   $\psi_{\pm}(x, k)$ as $|x|\to\infty$. For definiteness, we choose
  $\psi_-(x, k)= :\psi (x, k)$. We proceed from the Lippman--Schwinger equation \e{eq:WF2}. There is a significant difference between the cases $n\leq 2$ and $n>2$. Consider the resolvent kernel  \e{eq:R} for $z= k^n +i0$.
    If $n\leq 2$, then it  consists of oscillating terms, but,    for $n>2$,  additional  decaying terms  appear. We have to distinguish  oscillating   and decaying     parts in  the  right-hand side of \e{eq:WF2}.

    The decaying part $\psi_{\rm dec} (x,k)$  is defined by  the formula
    \begin{align}
\psi_{\rm dec} (x,k) = i n ^{-1}\sum_{\Im \kappa_{j}>0} \kappa_{j}^{-n+1}  \int_{-\infty}^x
 e^{i \kappa_{j} (x- y)} w(y,k) dy
\nonumber \\
 -i n ^{-1}\sum_{\Im \kappa_{j}<0} \kappa_{j}^{-n+1}  \int_x ^\infty 
 e^{i \kappa_{j} (y -x)}   w(y,k) dy
\label{eq:LS3}\end{align}
where $\kappa_{j}^n = k^n$ and 
    \begin{equation}
w (x,k) = ( V \psi(k) )(x).
\label{eq:LS4}\end{equation}
It follows from estimates \e{eq:WFDD} that
  \begin{equation}
| w(x,k)| \leq C (k) (1+|x|)^{-\rho} , \q \rho>1, \q \lambda= k^n\in \Omega .
\label{eq:WFDw}\end{equation}

 \begin{lemma}\label{dec}
 Function \e{eq:LS3} satisfies the condition
 \[
\lim_{|x|\to\infty}| \psi_{\rm dec} (x,k) |=0.
\]
  \end{lemma}

 \begin{pf}
  Consider, for example, one of the integrals in the first sum in \e{eq:LS3}. Obviously, it
 tends to $0$ as $x\to-\infty$ because  $w\in L^{1}({\Bbb R})$ and $e^{ \Im \kappa_{j} (y-x)}\leq 1$. If $x\to+\infty$, then we use the estimate  
  \begin{equation}
 \big|  \int_{-\infty}^x
 e^{i \kappa_{j} (x- y)}   w(y,k) dy\big| \leq e^{- \Im \kappa_{j} x/2} \int_{-\infty}^{x/2}
 |  w(y,k) |dy+  \int^{x}_{x/2} |  w(y,k) |dy, \q \Im \kappa_{j}>0, 
\label{eq:xd}\end{equation}
  where both terms on the right tend to zero.
   \end{pf}  
  
To define the oscillating  part $\psi_{\rm osc} (x,k) $ of $\psi (x,k) $, we introduce
the function
\begin{equation}
 r (x,k) = 1-  i n^{-1}  k^{-n+1}  \int_{-\infty}^x e^{-ik y} w(y,k) dy
\label{eq:OSC3}\end{equation}
if $n$ is odd.  If $n$ is even, we introduce two functions 
  \begin{equation}
  \begin{split}
r_{+}(x,k) &=1-  i n^{-1}  k^{-n+1}   \int_{-\infty}^x e^{-ik y} w(y,k) dy, \q k>0,
 \\
r_{+}(x,k) &=1 + i n^{-1}  k^{-n+1}   \int_{x}^\infty e^{-ik y} w(y,k) dy, \q k<0,
\end{split} 
\label{eq:ev1}\end{equation}
and
 \begin{equation}\begin{split}
r_{-}(x,k)& =-  i n^{-1}  k^{-n+1}   \int_{x}^\infty e^{ik y} w(y,k) dy, \q k>0,
 \\
r_{-}(x,k) &=  i n^{-1}  k^{-n+1}   \int_{-\infty}^x e^{ik y} w(y,k) dy, \q k<0 .
\end{split} \label{eq:ev3}\end{equation}
 In all these formulas $w$ is function \e{eq:LS4}.
Now we set
\begin{equation}
\psi_{\rm osc}  (x,k) =   e^{i k x} r (x,k) \q {\rm for}\; {\rm odd}\q n,
\label{eq:OSC1}\end{equation}
\begin{equation}
\psi_{\rm osc}  (x,k) =   e^{i k x} r_{+} (x,k)+ e^{-i k x} r_{-} (x,k) \q {\rm for}\; {\rm even}\q n.
\label{eq:OSC2}\end{equation}
With these definitions,  it follows from formulas \e{eq:R} and the Lippman--Schwinger equation \e{eq:WF2} that
   \begin{equation}
 \psi(x,k) =\psi_{\rm osc}  (x,k)  +\psi_{\rm dec} (x,k).
\label{eq:LS}\end{equation}

Since $w\in L^1 ({\Bbb R})$, the functions $ r  (x,k)$ and $ r_{\pm}  (x,k)$ have finite limits as $x\to\infty$ and as $x\to-\infty$.
Moreover,   we  have
$ r  (-\infty,k)=1$
for $n$ odd and
 % \begin{equation} 
 \[
    \begin{split}
r_{+}(-\infty,k) =1, \q  r_{-}(+\infty,k) =0\q {\rm if}\q k>0,
\\
r_{+}(+\infty,k) =1, \q  r_{-}(-\infty,k) =0\q {\rm if}\q k<0 ,
  \end{split}
  \]
%\label{eq:inf+}\end{equation}
for $n$ even. Let as usual $\lambda=k^n$.
We set
\begin{equation}
s(\lambda)= r  (+\infty,k) 
 \label{eq:DES}\end{equation}
for odd $n$
 and
\begin{equation}
   \begin{split}
s_{11}(\lambda)=r_{+} (+\infty,k), \q s_{21}(\lambda)=r_{-} (-\infty,k), \q k>0,
 \\
s_{12}(\lambda)=r_- (+\infty,k),\q s_{22}(\lambda)=r_{+} (-\infty,k), \q k<0,
   \end{split}
 \label{eq:SMX}\end{equation}
 for even $n$.

According to Lemma~\ref{dec} the function $\psi_{\rm dec} (x,k)$ does not contribute to  the asymptotics of $ \psi(x,k) $.
Therefore  representations  \e{eq:OSC1}, \e{eq:OSC2} and \e{eq:LS} yield  the following result.

 \begin{theorem}\label{SpThEF}
 Let estimates \e{eq:SR}  be satisfied.
  Then the functions $\psi (x, k) = \psi_{-}(x, k)$  defined by formula \e{eq:WF} satisfy the differential equation \e{eq:DE} 
  and have the following asymptotic behavior as $|x|\to \infty$:
  
  \begin{enumerate}[\rm(i)]
\item
If    $n$ is odd, then
\begin{equation}\left\{\begin{array}{lcl} 
\psi (x,k) =  e^{  ik x}+o(1) ,\quad x\ri -\infty,
\\  
\psi (x, k) =  s (\lambda)   e^{  ik x}+o(1) ,\quad x\ri \infty,
\end{array}\right.
\label{eq:WFSodd2}\end{equation}

 \item
 If    $n$ is even, then
 \begin{equation}\left\{\begin{array}{lcl} 
\psi(x,k) = s_{11}(\lambda) e^{   ik x}+o(1) ,\quad x\ri \infty,
\\
   \psi (x,k ) =e^{  ik x}+s_{21}(\lambda) e^{  -ik x}+o(1) ,\quad x\ri -\infty,
\end{array}\right.
\label{eq:WFS1}\end{equation}
for $k>0$ and
\begin{equation}\left\{\begin{array}{lcl} 
\psi (x,k) =e^{  ik x}+s_{12}(\lambda) e^{ -  ik x}+o(1) ,\quad x\ri \infty,
\\  
\psi (x, k) = s_{22}(\lambda) e^{ ik x}+o(1) ,\quad x\ri -\infty,
\end{array}\right.
\label{eq:WFS2}\end{equation} 
for $k<0$. 
 \end{enumerate}
   \end{theorem}
   
   Note that,  for odd $n$, expression  \e{eq:OSC3},  \e{eq:DES} for $s (\lambda)$ is the standard stationary representation of the scattering matrix. Therefore the asymptotic coefficient $s (\lambda)$ in \e{eq:WFSodd2} coincides with  $s (\lambda)$ defined  by \e{eq:S1odd}. Similarly, for even $n$, the   stationary representation of the scattering matrix is given by formulas \e{eq:ev1}, \e{eq:ev3} and \e{eq:SMX}. Therefore
      the asymptotic coefficients $s_{j \ell} (\lambda)$ in \e{eq:WFS1}   and \e{eq:WFS2} coincide with the entries of the scattering matrix \e{eq:SM}. 
   
 % these coefficients coincide with the entries of the scattering matrix \e{eq:SM}. The coefficient \e{eq:DES} coincides with $s(\lambda)$ defined by \e{eq:S1odd}.
   
%   \begin{equation}\left\{\begin{array}{lcl} 
%\psi(x,k) =   e^{   ik x}  -  i n^{-1 } |k|^{-n+1 }e^{  - i|k| x} \int_{-\infty}^\infty e^{  i|k| y}(V \psi(k)) (y) dy+o(1) ,\quad x\ri -\infty,
%\\   \psi (x,k ) =e^{  ik x}-  i n^{-1 } |k|^{-n+1 }e^{   i|k| x} \int_{-\infty}^\infty e^{ -  i|k| y}(V \psi(k)) (y) dy +o(1) ,\quad x\ri +\infty. \end{array}\right. \label{eq:WFSS}\end{equation}
   
 %We have defined the asymptotic coefficients $s_{j,\ell} (\lambda)$ by the integral representations  \e{eq:ev1} -- \e{eq:ev4}. It can be standardly shown that these coefficients coincide with the entries of the scattering matrix \e{eq:SM}. The coefficient \e{eq:DES} coincides with $s(\lambda)$ defined by \e{eq:S1odd}.

 Similarly to the operators of second order,  the solution $\psi(x, k) $   of differential equation  \e{eq:DE} describes the plane wave coming from minus infinity if $k>0$ (from plus infinity if $k<0$). If $n$ is even, then the numbers  $s_{11}(\lambda)$ and $s_{21}(\lambda)$ in \e{eq:WFS1} (the numbers  $s_{22}(\lambda)$ and $s_{12}(\lambda)$ in \e{eq:WFS2}) are interpreted as  the corresponding transmission and reflection coefficients. If $n$ is odd, then according to   \e{eq:WFSodd2} the reflection coefficients are zeros.

Putting together Theorem~\ref{SpThEF} with formula \e{eq:SMev}, we find also the asymptotics of the functions $\psi_{+}(x,k)$ as $x\to\infty$ and as $x\to-\infty$.

 \medskip

   {\bf 2.5.} 
 %  In view of our applications to Hankel operators, 
      Let us specially discuss a particular case when   the coefficient
    $b_{n-1} (x)$ is long-range while other coefficients
    $b_{1} (x), \ldots, b_{n-2} (x)$ are short-range.  Recall that, by a gauge transformation,  one can  get rid of the term $b_{n-1}(x)D^{n-1}$ in \e{eq:VB}. Indeed,
 let
 \begin{equation}
 ( {\cal J}  f)(x) = e^{i \beta (x)} f(x)
\label{eq:Ga}\end{equation}
where
\begin{equation}
\beta (x)= - \frac{1}{n}\int_{0}^x b_{n-1}(y)dy.
\label{eq:Ga1}\end{equation}
Since $D{\cal J} ={\cal J} (D + \beta')$, the operator 
  \begin{equation}
 \wt{B} =  {\cal J} ^* B  {\cal J}   
  \label{eq:Ga2}\end{equation}
has again the form  \e{eq:VB} with   $\wt{b}_{ n-1}(x)= 0$. Moreover, the coefficients     
  $\wt{b}_{m}(x)$, $m=0,\ldots,  n-2$, of the operator $\wt{B}$ are short range if this is true for  $ {b}_{n-1}^2(x)$ and the derivatives of $ {b}_{n-1}(x)$. In this case   we can directly apply the results of subs.~2.2,  2.3
 and 2.4 to the operator $\wt{B}$   and then    carry over these results   to the operator $B$. Indeed, it follows from  definition \e{eq:WLR} and relation \e{eq:Ga2} that the wave operators
   \begin{equation}
  W_\pm :=  W_\pm(B , B_{0};   {\cal J} ) = {\cal J} W_\pm(\wt{B} , B_{0}  ) , \q B_{0} =D^n,
  \label{eq:WWS}\end{equation}
  exist.   For the corresponding scattering operators,  we have the identity
    \begin{equation}
 {\cal S}:= W_+ (B , B_{0};   {\cal J} ) ^*  W_- (B , B_{0};   {\cal J} )  =   W_+(\wt{B} , B_{0} )^* W_-(\wt{B} , B_{0} ) .
  \label{eq:WWS1}\end{equation}
Thus the scattering matrices for the pair $B_{0}$, $\wt{B}$ and the triple $B_{0}$, $B $, ${\cal J}$  are the same. 
  
In view of   \e{eq:WF} the eigenfunctions of the operator $\wt{B}$ are defined by the equation
\begin{equation}
\wt\psi_{\pm } (k)   =\psi_{0}  (k)  - \wt{R} (k^{n}\mp i0)(\wt{B}-B_{0} ) \psi_{0}  (k)  , \q \lambda= k^n\in\Omega,
\label{eq:WFt}\end{equation} 
where $\wt{R} (z)=(\wt{B}-z)^{-1}$. According to   \e{eq:Ga}, \e{eq:Ga2} it is natural to define
  eigenfunctions of the operator $B $   by the relation 
  \begin{equation}
  \psi_{\pm} (x, k)=  e^{i \beta(x)}\wt{\psi}_{\pm} (x, k).
  \label{eq:Ga3}\end{equation}

This leads to the following result.

      \begin{proposition}\label{BL}
      Let the operator $B$ be given by formula \e{eq:VB}.       Assume that 
    estimates  \e{eq:SR} with some $\rho>1$   are true  for $m=0, \ldots, n-2$.   With respect to the coefficient $b_{n-1} (x)$,  we assume that
    $b_{n-1} \in C^{n-1}$, $b_{n-1}  (x)=O ( |x|^{- \rho/2} )$    and 
    $b_{n-1}^{(p)} (x)=O ( |x|^{-\rho })$ 
  for   $p=1,\ldots, n-1$.  Then:
\begin{enumerate}[\rm(i)] 
  \item
  All conclusions of Theorem~\ref{SpTh} for the operator $B$ are true.
 \item
 If $\cal J$ is defined by \e{eq:Ga}, \e{eq:Ga1}, then    the  wave operators
$  W_\pm(B , B_{0};   {\cal J} ) $ exist. Moreover, if $\wt{B}$ is defined by  \e{eq:Ga2}, then the scattering matrices for the pair $B_{0}$, $\wt{B}$    and the triple $B_{0}$, $B $, $\cal J$  are the same.
 \item
   The functions  $  \psi_{\pm} (x, k)$ defined by \e{eq:Ga3} satisfy differential equation   \e{eq:DE} and their asymptotic behavior as $|x| \to \infty$ is determined by formulas  \e{eq:WFSodd2}  $($for $n$ odd$)$ and \e{eq:WFS1}, \e{eq:WFS2}  $($for $n$ even$)$  
  where the right-hand sides acquire the additional factor $e^{i \beta (x)}$.  Relation  \e{eq:est2B}  remains also true.
\end{enumerate} 
    \end{proposition}

 \medskip

   {\bf 2.6.} 
   In view of our applications to Hankel operators, 
   let us finally discuss the case of the coefficients satisfying the condition
 \begin{equation}
b_{m}(x)= \ov{b_{m}(-x)}, \q m=0,\ldots, n-1.
\label{eq:Compl}\end{equation}
Then the operator $B$ commutes 
with the antilinear  involution ${\sf C}$  defined by the relation
 \begin{equation}
 ( {\sf C} f)(x)= \ov{f(- x)}.
\label{eq:Complx}\end{equation}
In this case the integral kernel of the resolvent $R(z)= (B-z)^{-1}$ satisfies the identity $R(x,y;z)= R(-y,-x;z)$.
 Let the operators $\cal J$, $W_{\pm}$ and $\cal S$ be defined by equalities  \e{eq:Ga},  \e{eq:WWS} and \e{eq:WWS1}, respectively. Since ${\sf C} \cal J =\cal J {\sf C} $, ${\sf C} {\cal J}^* = {\cal J}^* {\sf C}$ and  $ {\sf C} e^{iBT}=e^{- iBT} {\sf C}$, we see that ${\sf C} W_{\pm} = W_{\pm}  {\sf C}$, ${\sf C} W_{\pm}^* = W_{\pm}^* {\sf C}$
 and therefore
 \begin{equation}
{\sf C} \cal S = { \cal S}^* {\sf C} . 
\label{eq:Comply}\end{equation}
Obviously, $\Phi {\sf C} \Phi^*= {\cal C}$ is the complex conjugation, that is, $({\cal C}  g)(k)=\ov{g(k)}$. Therefore the identity \e{eq:Comply} implies that ${\cal C} \wh{\cal S}  {\cal C}  =\wh{\cal S}^*$. In terms of the scattering matrices this result can be reformulated in the following form.

 % Recall that, for $n$ even,  the scattering matrix $S(\lambda)$ was defined  by relations \e{eq:S1odd} and \e{eq:S1}.
%  Since $ \cal C\Phi f= \Phi\bar{f}$, we see that
%  $$  \ov{Y\Phi{ \cal S}\Phi^{*} Y^{*}f}=Y\Phi { \cal C}{ \cal S}{ \cal C} \Phi^{*} Y^{*}\bar{f}.  $$
%  Therefore     the following result is a direct consequence of   \e{eq:Comply}. 

    \begin{proposition}\label{BL1}
    Let the conditions of Theorem~\ref{SpTh} or of Proposition~\ref{BL} be satisfied.
    Under assumption  \e{eq:Compl}, we have the identity
    \begin{equation}
   \ov{S(\lambda)}= S(\lambda)^*.
\label{eq:Complz}\end{equation}
    \end{proposition}
    
 For $n$    odd, the scattering matrix  is a scalar function $s(\lambda)$ so that
        this identity only means that $| s(\lambda)|=1$, and hence     it
      does not impose any additional  restrictions on $s(\lambda)$.
      On the contrary,    for $n$ even, the identity \e{eq:Complz} implies that $s_{12} (\lambda)= s_{21} (\lambda)$, that is, the right and left reflection coefficients coincide. This result  may be compared with the fact that for the Schr\"odinger operator with a real potential, the transmission coefficients coincide, that is,  $s_{11} (\lambda)= s_{22} (\lambda)$ (the paper \cite{F} by L.~D.~Faddeed or  the book \cite{Ya}, \S 5.1).

    \section{Long-range perturbations}
    
  % \medskip
 
     {\bf 3.1.} 
   Let us now consider the operator $B$ defined by formula \e{eq:VB} where the coefficients $b_m  (x)$ are long-range.
   To be precise, we suppose  that 
   \begin{equation}
| b_m^{(p)} (x)| \leq C_{p} (1+|x|)^{-\rho -p}, \q \rho >0,\q m=0,1, \ldots, n-1,
 \label{eq:LR}\end{equation}
 for all $p=0,1,\ldots$ although this condition for some finite number (depending on $\rho$) of $p$ is sufficient. The number $\rho$ may be arbitrary small.  Thus  the functions $b_m (x) $ decay at infinity but perhaps very slowly. 
 Nevertheless we have the following result.
 
  \begin{theorem}\label{SpThlong}
  Under assumption \e{eq:LR} all conclusions $(i), (ii), (iii)$ and $(iv)$ of
  Theorem~\ref{SpTh} are true.  
  \end{theorem}
 
 In the long-range case   the method exposed in subs.~2.1 and 2.2 and relying on the resolvent equation \e{eq:res}  is no longer applicable. However, similarly to second order differential operators $B$, 
  assertions (i) and (ii) can be proven  by  using the Mourre method \cite{Mo1}. 
  
We recall it here very briefly.  Let   ${\Bbb D}= xD + D x$, and let as before $B_{0}= D^{n}$, $V=B-B_{0}$. The Mourre method \cite{Mo1} relies on the commutation relation
    \[
    i [B_{0}, {\Bbb D}]: = i B_{0} {\Bbb D} -i  {\Bbb D} B_{0}=  2n  B_{0}.
    \]
       The commutator $[ V, {\Bbb D}]$ is a differential operator of order $n-1$ with coefficients decaying at infinity.  This allows one to estimate the commutator  $i [ B, {\Bbb D}]$ from below.  Let $E(X) $ be the spectral projector of the operator $B$ corresponding to a set $X\subset {\Bbb R}$. If $\lambda\neq 0$ is not an eigenvalue of the operator $B$, then for a sufficiently small $\delta>0$ and $X=(\lambda- \delta,\lambda+\delta)$
        the Mourre estimate 
     \[
    i E (X)[B, {\Bbb D}] E (X)\geq c E (X), \q c=c(\lambda) >0,
  \]
  %\label{eq:Mourre }\end{equation}
  holds. 
  
  To deduce from this estimate the   assertions (i) and (ii) of
  Theorem~\ref{SpThlong}, one needs to consider the second commutator
  \[
 [ [B, {\Bbb D}] ,{\Bbb D}]=-4 n^2 B_{0} + [ [V, {\Bbb D}] ,{\Bbb D}] 
 \]
 where again $ [ [V, {\Bbb D}] ,{\Bbb D}] $ is a differential operator of order $n-1$ with coefficients decaying at infinity. 
 It follows that the operator $ [ [B, {\Bbb D}] ,{\Bbb D}] (B_{0}+i)^{-1}$ is bounded. As shown by E.~Mourre, this fact is sufficient for the proof of  the   assertions (i) and (ii). Its presentation can be found in various books; see, e.g., Chapter~4 of  \cite{CFKS}, Chapter~7 of  \cite{ABG} or \S 6.9 of \cite{Ya}.

    We note that the Mourre method applies also in the short-range case.

 \medskip
 
     {\bf 3.2.} 
     Similarly to the short-range case,   the proof of assertion (iii) can be obtained by tools of the scattering theory. However for long-range coefficients  $b_{m}(x)$, the operator $B$ cannot be considered as a perturbation of $B_{0}$. In particular, this means that
    the approximation as $T\to\infty$ of  $\exp(-iB T)$   by the free evolution $\exp(-iB_{0}T)$  is not sufficient or, to put it differently,   the usual wave operators $W_{\pm} (B,B_{0})$ for the pair $B_{0}$, $B$ do not exist. Therefore  we have to replace them by   more 
general wave operators \e{eq:WLR}
where the identification $J$ takes into account
 the behavior of the coefficients $b_{m} (x)$ of $B$ as $|x|\to\infty$.  We choose $J$ as a  pseudodifferential operator whose symbol depends on the coefficients of the operator $B$. We follow  here a very simplified version of the scheme suggested in \cite{Y10} for the multidimensional Schr\"odinger operator.
 
 %   This   is discussed in the Appendix.

%  Now we assume condition  \e{eq:LR} with some $\rho\in (0,1]$ to be satisfied.

 We seek $J$ as a pseudodifferential operator
 \begin{equation}
( J f) (x)= (2\pi)^{-1/2} \int_{-\infty}^\infty e^{i \phi(x,k)} \chi (k^n) \hat{f} (k) dk
\label{eq:WLR1}\end{equation}
where $\chi \in C_{0}^\infty ({\Bbb R}\setminus
\{0\})$,
 \begin{equation}
\phi(x,k)= xk+ \vartheta(x,k)
\label{eq:WLR2}\end{equation}
and $\vartheta(x,k)=o(x)$ as
  $|x|\to\infty$. The function $\vartheta (x,k)$ is  introduced to handle   the slow decay of the coefficients $b_{m} (x)$ at infinity. The auxiliary function   $\chi $ 
allows us to neglect high and low energies. Now the ``perturbation''  ${\bf V}:=B J-J B_{0}$ is given by the formula
 \begin{equation}
( {\bf V} f) (x)= (2\pi)^{-1/2} \int_{-\infty}^\infty e^{i  x k} e^{i \vartheta(x,k)} {\bf v}(x,k)\chi  (k^n) \hat{f} (k) dk
\label{eq:WLR3}\end{equation}
where
 \begin{equation}
{\bf v} (x,k)= e^{-i \phi(x,k)} (D^n+ \sum_{m=0}^{n-1}b_{m}(x) D^m -k^n)  e^{i \phi(x,k)}. 
\label{eq:WLR4}\end{equation}
We have to construct a function $\vartheta(x,k)$ such that  
 \begin{equation}
  \vartheta^{(l)}(x,k)= O (|x|^{1-l- \rho}),  \q  l=0,1,\ldots,
\label{eq:WLR9aa}\end{equation}
and  
 \begin{equation}
 {\bf v}^{(l)}(x,k)= O (|x|^{-1-l- \rho}),  \q  l=0,1,\ldots,
\label{eq:WLR5x}\end{equation}
as $|x|\to\infty$ uniformly in $k$ on compact subsets of $ {\Bbb R}\setminus\{0\}$. We also require that the same estimates
\e{eq:WLR9aa} and \e{eq:WLR5x}  are true for all derivatives  of $ \vartheta^{(l)}(x,k)$ and $ {\bf v}^{(l)}(x,k)$ in $k$. A function $ {\bf v}(x,k)$ satisfying these conditions will be called short-range. Condition \e{eq:WLR9aa} means that the  function \e{eq:WLR2} is close to the linear function $xk$ as $|x|\to\infty$.

Suppose that $\vartheta(x,k)$ satisfies condition \e{eq:WLR9aa}, and 
let $\phi(x,k)$ be defined by \e{eq:WLR2}. Then 
for all $m$ we have
 \[
e^{-i \phi(x,k)} D^m e^{i \phi(x,k)}=   \big( k+ \sigma(x,k)\big)^m+\tau_{m} (x,k)
\]
%\label{eq:WLR6}\end{equation}
where $  \sigma(x,k)= \vartheta'(x,k)$ and $\tau_{m} (x,k)$
is a short-range function.    Substituting this expression into \e{eq:WLR4}, we see that, up to short-range terms,
 \begin{equation}
{\bf v} (x,k)=n k^{n-1}  \sigma(x,k) + \sum_{m=2}^{n}\tbinom{n}{m}   \sigma(x,k)^m k^{n-m} +
 \sum_{m=0}^{n-1}b_{m}(x) (k+\sigma(x,k))^m. 
\label{eq:WLR7}\end{equation} 
Let us construct an approximate solution of the equation ${\bf v} (x,k)=0$ for $\sigma (x,k)$ by iterations. We set  $\sigma_{0}(x,k)=0$,
 \begin{equation}
 n k^{n-1}  \sigma_{1}(x,k)   =-
 \sum_{m=0}^{n-1}b_{m}(x) k^m
\label{eq:WLp8}\end{equation}
  and 
 \begin{equation}
 n k^{n-1} \sigma_{j+1}(x,k) =- \sum_{p=2}^{n}\tbinom{n}{p}  \sigma_{j}(x,k)^p k^{n-p}  -
 \sum_{m=0}^{n-1}b_{m}(x) (k+ \sigma_{j}(x,k))^m 
\label{eq:WLR8}\end{equation}
for all $j\geq 1$.  Let $\sigma (x,k)=\sigma_{j}(x,k)$ in \e{eq:WLR7}. Then
  \begin{equation}
{\bf v} (x,k)=n k^{n-1} \big( \sigma_{j}(x,k) - \sigma_{j+1}(x,k) \big).
\label{eq:bfv}\end{equation}

We will  check that
  \begin{equation}
 \sigma_{j }^{(l)}(x,k) -  \sigma_{j-1}^{(l)}(x,k)= O (|x|^{-l-  j \rho}) 
\label{eq:WLR9}\end{equation}
and hence
  \begin{equation}
  \sigma_{j}^{(l)}(x,k)= O (|x|^{-l- \rho})
\label{eq:WLR9a}\end{equation}
for all $j =1, 2, \ldots$ and $l=0,1,\ldots$. For $j=1$, \e{eq:WLR9} is a direct consequence of definition \e{eq:WLp8} and assumption \e{eq:LR}.
Let us justify the passage from $j$ to $j+1$.
It follows from \e{eq:WLR8} that
 \begin{multline}
 n k^{n-1} ( \sigma_{j+1}(x,k) -\sigma_{j}(x,k) )=- \sum_{p=2}^{n}\tbinom{n}{p}   \big(\sigma_{j}(x,k)^p   - \sigma_{j-1}(x,k)^p \big) k^{n-p}
 \\
- \sum_{m=0}^{n-1}b_{m}(x)\big( (k+ \sigma_j(x,k))^m-  (k+\sigma_{j-1}(x,k))^m\big) 
\label{eq:WLQ}\end{multline}
where
\begin{equation}
\big(k+ \sigma_j(x,k)\big)^m-  \big(k+ \sigma_{j-1}(x,k)\big)^m=  \sum_{p=1}^{m}\tbinom{m}{p}  
\big( \sigma_j(x,k)^p    - \sigma_{j-1}(x,k)^p\big) k^{m-p}. 
\label{eq:WLQ2}\end{equation}
Using \e{eq:WLR9} and \e{eq:WLR9a} we see that for all $p$
 \[
  D^l \big(\sigma_{j}(x,k)^p - \sigma_{j-1}(x,k)^p \big)
= O (|x|^{-l -(j + p -1)\rho}) 
\]
%\label{eq:WLQ1}\end{equation}
 % In view of \e{eq:WLQ1} this expression is  $O (|x|^{ -j  \rho})$, and its derivative of order $l$ is $O (|x|^{-l -j  \rho})$.
 Substituting these estimates into the right-hand sides of \e{eq:WLQ} and  \e{eq:WLQ2} and taking into account condition  \e{eq:LR}  on the coefficients $b_{m}(x)$, we get estimate \e{eq:WLR9} with $j$ replaced by $j+1$.

The following result is a direct consequence of estimates  \e{eq:WLR9}, \e{eq:WLR9a} and representation \e{eq:bfv} for ${\bf v} (x,k)$.

\begin{lemma}\label{LRSt}
Let the assumption  \e{eq:LR} with $\rho\in (0,1]$ hold. 
  Set $\sigma(x,k)= \sigma_{j}(x,k)$ where $j > \rho^{-1}$. Then the function
  \[
  \vartheta (x,k)=\int_{0}^x  \sigma(y,k) dy.
  \]
satisfies   condition \e{eq:WLR9aa}. Let the functions $\phi(x,k)$ and  ${\bf v} (x,k)$ be defined by formulas 
\e{eq:WLR2} and  \e{eq:WLR4}. Then estimates \e{eq:WLR5x} hold. 
  \end{lemma}

  \medskip 
 
{\bf 3.3.}
 The operator $J$ defined by formulas \e{eq:WLR1},  \e{eq:WLR2} is a pseudodifferential operator with symbol
 \[
 {\bf j}(x,k)= e^{i \vartheta(x,k)}    \chi  (k^n).
 \]
 Here and below $\vartheta(x,k)$ is the function constructed in Lemma~\ref{LRSt}.
 
 Recall that the H\"ormander class ${\cal S}^m_{\rho, \d}$, $\rho\in [0,1], \d \in [0,1]$, of symbols $ {\bf a}(x,k)$ compactly supported in $k$ is distinguished by estimates
  \begin{equation}
|(\partial_{x}^\alpha\partial_{k}^\beta{\bf a})(x,k)| \leq C_{\alpha,\beta} (1+|x|)^{m-\alpha\rho+\beta \d}, \q \alpha,\beta=0,1, \ldots.
\label{eq:PD}\end{equation}
Pseudodifferential calculus  (see the books \cite{Hor3, Sh}) works smoothly in this class if $\rho>1/2>\d$. Otherwise, a specific study is required. It follows from estimates \e{eq:WLR9aa} that
 ${\bf j}\in {\cal S}^0_{\rho, 1-\rho}$. Therefore for $\rho>1/2$, the next result is a consequence of general results on pseudodifferential  operators. On the contrary, for $\rho\leq 1/2$ the oscillating nature of symbol \e{eq:PD} has to be taken into account (see Theorems~2.12 and 2.13  in  \cite{Y13}).
 
 \begin{lemma}\label{LRStj}
Under the assumption  \e{eq:LR}   the operator $J$ defined by formulas \e{eq:WLR1},  \e{eq:WLR2} is bounded. The operators $J^* J - \chi (B_{0})^2$ and $J J^*   - \chi (B_{0})^2$ are compact.
  \end{lemma}

Now we are in a position to prove the existence of the wave operators
 \e{eq:WLR} and \e{eq:WLRad}.  Relation  \e{eq:WLR3} means that the perturbation $ {\bf V}$
  is a pseudodifferential operator with symbol
 \[
 {\sf v}(x,k)= e^{i \vartheta(x,k)} {\bf v}(x,k)\chi  (k^n).
 \]
 According to  Lemma~\ref{LRSt} we have $ {\sf v} \in {\cal S}^{-1-\rho}_{\rho, 1-\rho}$ which directly implies that
 the operator $\la x\ra^{\sigma} (B J-J B_{0})\la x\ra^{\sigma}$   is bounded if   $\sigma=(1+\rho)/2$.
 By part (ii) of Theorem~\ref{SpThlong} the operator $\la x\ra^{-\sigma}$  is $B_{0}$- and $B$- smooth in the sense of Definition~\ref{Kato}. Therefore the next assertion is a  consequence of Proposition~\ref{smooth}. 
 
   \begin{theorem}\label{LRSt1}
Let $B_{0}=D^n$, let $B$ be  the operator  \e{eq:VB}, and let the operator $J$ be defined by formulas \e{eq:WLR1},  \e{eq:WLR2}.  Then under the assumption  \e{eq:LR}  the wave operators \e{eq:WLR} and \e{eq:WLRad} exist. 
  \end{theorem} 
 
Put $W_{\pm} = W_{\pm} (B,B_{0}; J)$.  Lemma~\ref{LRStj} ensures that 
 \begin{equation}
W_{\pm}^* W_{\pm} = \chi (B_{0})^2 \q {\rm and}\q  W_{\pm} W_{\pm}^*  = \chi  (B)^2 .
\label{eq:PDWO}\end{equation}
Recall that the open set $\Omega$ was defined by relations \e{eq:omega}.
 Let  $X\subset \Omega$ be a compact interval, and let $E_{0} (X)$, $E (X)$ be the spectral projectors of the operators $B_{0}, B$. Choose a function $\chi  \in C_{0}^\infty ({\Bbb R}\setminus\{0\})$ such that $\chi  (\lambda)=1$ for $\lambda\in X$. Then it follows from \e{eq:PDWO}
 that for all $f\in L^2 ({\Bbb R})$,
  \[
\|W_{\pm} E_{0}(X)f\|= \| E_{0}(X)f\| \q {\rm and}\q  \|W_{\pm}^* E (X)  f\|= \| E (X)   f\|  .
\]
%\label{eq:PDWO1}\end{equation}
Using that $W_{\pm} B_{0}=B W_{\pm}$, we see that the operators $B_{0}E_{0}(X)$ and $B E(X)$ are unitary equivalent. Since $X$ is arbitrary, this implies the unitary equivalence of the operators $B_{0} $ and $B^{(ac)} $ and concludes the proof of statement (iii) of Theorem~\ref{SpThlong}.

The proof of statement (iv) is practically the same as in the short-range case. The only difference is that we now have to replace   the phase $k_{j}x$ in \e{eq:EigF} by a more general phase function with the same asymptotic behavior at infinity.

  \medskip 
 
{\bf 3.4}.
Theorem~\ref{LRSt1} allows us to find the asymptotics of $\exp(-i B T) f$ as $T\to\pm\infty$  for all $f\in{\cal H}^{(ac)}$.
Let  again $X\subset \Omega$ be a compact interval and $f=E(X)f$. Since
\begin{equation}
f=  W_{\pm} (B , B_{0};  J) f_{ 0}^{(\pm)} \q  {\rm where}\q f_{ 0}^{(\pm)}= W_{\pm} (B_{0}, B; J^*) f,
\label{eq:StPxb}\end{equation}
it follows from definition  \e{eq:WLR} that
 \begin{equation}
\lim_{T\to\pm\infty}\| \exp(-i B T) f -J  \exp(-i B_{0}T) f_{ 0}^{(\pm)}\|=0 .
\label{eq:StPx}\end{equation}
According to  \e{eq:WLR1} we have
\begin{equation}
( J  \exp(-i B_{0} T) f_{ 0}^{(\pm)} ) (x)= (2\pi)^{-1/2} \int_{-\infty}^\infty e^{i \phi(x,k)}e^{-i k^n T}  \chi (k^n) \hat{f}_{ 0}^{(\pm)} (k) dk 
\label{eq:WLR1b}\end{equation}
where $\phi(x,k)$ is given by equation \e{eq:WLR2}.

Let us use   the stationary phase method for $T\to\pm\infty$ (see, e.g.,  Lemma~\ref{StPh} where the role of $T$ is played by a large parameter $N$) with the phase function
$
\omega (x,k,T)= \phi(x,k)T^{-1}-k^n
$
depending on parameters $x$ and $T$. In view of \e{eq:WLR2} the stationary points of integral \e{eq:WLR1b} are determined by the equation
 \begin{equation}
x T^{-1}-n k^{n-1}   + T^{-1}\partial\vartheta (x,k)/\partial k=0
\label{eq:StPh3}\end{equation}
for $k = k(x,T)$.
Condition \e{eq:WLR9aa} on derivatives of $\vartheta (x,k)$ in $k$ allows us to obtain
the solutions of   equation \e{eq:StPh3}   by successive approximations starting from $k_{0}$ such that $n k_{0}^{n-1} T=x$. Note also that
 \begin{equation}
\partial^2 \omega (x,k, T)/\partial k^2 = -n (n-1) k^{n-2} + T^{-1}\partial^2\vartheta (x,k)/\partial k^2.
\label{eq:StPhT}\end{equation}

If $n$ is even, then   equation \e{eq:StPh3} has the unique solution $k (x,T)$ for all $x$ and sufficiently large $T$ and
 \begin{equation}
k (x,T)= (x/ (nT))^{\frac{1}{n-1}} + O(|T|^{-\rho})
\label{eq:kk}\end{equation}
uniformly   in  $x/T$ in  compact sets of ${\Bbb R}\setminus\{ 0\}$.
  Thus applying formula \e{eq:BHs1} to integral \e{eq:WLR1b} and using \e{eq:StPhT} we find that
 \begin{multline}
(J \exp(-i B_{0} T) f_{ 0}^{(\pm)})(x)= (n(n-1))^{-1/2}  e^{\mp \pi i/4} |T |^{-1/2}e^{i\Xi (x, T)}
\\
\times | k (x,T) |^{-\frac{n-2}{2}}  \hat{f}_{ 0}^{(\pm)} (k (x,T))+ \epsilon^{(\pm)}(x, T)
\label{eq:StPh}\end{multline}
where     the phase  
  \begin{equation}
  \Xi (x, T)= xk(x, T)- k(x,T)^n T + \vartheta (x, k(x, T))
\label{eq:StPh2}\end{equation}
and $\|\epsilon^{(\pm)}(T)\|\to 0$ as $ T \to\pm\infty$. Formula \e{eq:StPh} can be simplified if one  takes an approximate solution   of   equation \e{eq:StPh3} for $k(x, T)$. According to \e{eq:kk}    $|k(x, T)|$ in \e{eq:StPh} can be replaced by $|x/(nT)|^{\frac{1}{n-1}}$ and the phase \e{eq:StPh2}  can be written as
 \begin{equation}
  \Xi (x, T) =  \Xi_{0} (x, T)   +  \Theta (x/(nT) ,T) 
\label{eq:StPh2w}\end{equation}
with 
 \begin{equation}
  \Xi_{0} (x, T) = (n-1)   \big| x/ (n T)  \big|^{\frac{n}{n-1}} T    
\label{eq:Xi}\end{equation}
and some function  $ \Theta$ satisfying the condition
   $\Theta (u,T) = O(|T|^{1-\rho})$ for $u$ in  compact sets of ${\Bbb R}\setminus\{ 0\}$. 
   
   Let us   introduce a unitary
  operator  $Y= Y_{n}$   in $L^2 ({\Bbb R})$ by the formula
 \begin{equation}
(Y f )(x)=  (n-1)^{-1/2} | x |^{-\frac{n-2}{2(n-1)}} \hat {f} (\sgn x  |x|^{\frac{1}{n-1}} ).
\label{eq:StPh1}\end{equation} 
Then it follows from   \e{eq:StPh} that
 \begin{equation}
(J \exp(-i B_{0} T) f_{ 0}^{(\pm)})(x)=  e^{\mp \pi i/4}   | n T|^{-1/2}e^{i\Xi (x, T)}    (Y f_{ 0}^{(\pm)}) (x/ (nT))+ \epsilon^{(\pm)}(x, T)
\label{eq:StPhX}\end{equation}
where again $\|\epsilon^{(\pm)}(T)\|\to 0$ as $ T \to\pm\infty$.
 Putting together relations  \e{eq:StPx} and  \e{eq:StPhX}, we obtain the asymptotics of $\exp(-iB T) f$ as $T\to\pm \infty$ for vectors $f$ such that $f=E (X)f$.
  The unitarity of $\exp(-iB T) f$ allows us to extend it  to all $f\in{\cal H}^{(ac)} $. 
  
   If $n >1$ is odd, then equation \e{eq:StPh3} has no solutions for   $x/ T< 0$ and it has two solutions  $k_{j} (x,T)$ close
   to  $ (-1)^{j} (x/(n T)) ^{1/(n-1)}$, $j=1, 2$, for   $x/ T>0$. In this case we have two stationary points $k_{j} (x,T)$, and  instead of \e{eq:StPh2w} we have two phase functions
    \begin{equation}
  \Xi_{j} (x, T) = (-1)^{j}  \Xi_{0} (x, T)  + \Theta_{j} (x/ (n T),T),  \q j=1, 2,
\label{eq:StPhdd}\end{equation}
where   again $ \Theta_{j} (u, T) = O(| T|^{1-\rho})$ for $u$ in  compact sets of ${\Bbb R}_{+}$. Otherwise the considerations are quite similar to the case of even $n$.
   
   Let us state the results obtained.
 
 \begin{theorem}\label{qqs}
Let   $B$ be  the operator  \e{eq:VB}, and let the assumption  \e{eq:LR}  be satisfied. Suppose that $f\in{\cal H}^{(ac)}$ and  $f_{ 0}^{(\pm)}$ is defined by \e{eq:StPxb}.   If $n$ is even, then  
 \begin{equation}
( \exp(-i B T) f )(x)= e^{\mp \pi i/4} | n T|^{-1/2}e^{i\Xi (x, T)} (Y f_{ 0}^{(\pm)})(x/ (  nT ))+ \epsilon^{(\pm)}(x,T)
\label{eq:StPhB}\end{equation}
where  $\|\epsilon^{(\pm)}(T)\|\to 0$ as $T\to\pm\infty$. If $n$ is odd,  $n >1$ and $T \to\pm\infty$, then
 \begin{equation}
( \exp(-i B T) f )(x)= | nT|^{-1/2}  \sum_{j= 1, 2}  e^{\mp (-1)^{j }\pi i/4}  e^{i\Xi_{j} (x, T)} (Y f_{ 0}^{(\pm)})( (-1)^{j}x/ (n T ))+ \epsilon^{(\pm)}(x, T),
\label{eq:StPhBod}\end{equation}
where $\pm x>0$  and $\|\epsilon^{(\pm)}(T)\|_{L^{2} ({\Bbb R}_{\pm})}\to 0$. Moreover, $\| \exp(-i B T) f \|_{L^{2} ({\Bbb R}_{\mp})}\to 0$ as $T\to\pm\infty$.
  \end{theorem}

Of course these asymptotic formulas remain true but can be significantly simplified under the assumptions of Section~2. If condition \e{eq:SR}  is satisfied, then $J=I$  and we can set $\Theta=\Theta_{1}=\Theta_{2} =0$ in \e{eq:StPh2w} and \e{eq:StPhdd}.  
  Under the assumptions of Proposition~\ref{BL} al these functions equal $\beta(x)$ defined by \e{eq:Ga1}; in this case $J={\cal J}$ is given by \e{eq:Ga}.

  % In the general case, \e{eq:StPh2w} gives the asymptotics of the phase  function $  \Xi (x,T)$ as $| T|\to\infty$ if  
% $x/ T$ belongs to   compact sets in ${\Bbb R}\setminus\{ 0\}$.  

   %************************************************************
\section{Degenarate differential operators}  
%***********************************************************

 {\bf 4.1.} 
  Let us now consider differential   operators $A$ defined by formula \e{eq:LOGz} where $Q(X)$ is polynomial \e{eq:LO4}
  with $q_{n}=1$. Thus we have
    \begin{equation}
 (A g) (\xi)= v(\xi) \sum_{m =0}^n q_m D_{\xi}   \big(v(\xi)  g(\xi)\big) ,\q D_{\xi} =-i d/d\xi, \q q_{n} =1.
 \label{eq:VV}\end{equation}
   We do not suppose that the function $v(\xi)$ is given by the formula 
 \e{eq:LOGZ1}.
Instead we accept very general

 \begin{assumption}\label{qq}
 The function  $v \in C^\infty ({\Bbb R})$,  $v(\xi)>0$ and $v(\xi)\to 0$   as $|\xi | \to \infty$.  
  \end{assumption}

Let us make the change of variables 
  \begin{equation}
x= x(\xi)= \int_0^\xi v(\eta)^{-2/n} d\eta
 \label{eq:V2}\end{equation}
  and define the operator $ \mathrm {L}$ by the formula  
     \begin{equation}
 (\mathrm {L} f) (\xi)=  x'(\xi)^{1/2} f( x(\xi))
  \label{eq:V1}\end{equation}
Under Assumption~\ref{qq} we have $x'(\xi)=v(\xi)^{-2/n}>0$ and $x (\xi) \to \pm\infty$ as $\xi\to \pm\infty$.
The operator $\mathrm {L} : L^2 (\Bbb R) \to L^2 (\Bbb R) $ is unitary.

 Let us calculate the operator
  \begin{equation}
  B=\mathrm {L}^* A \mathrm {L}.  
   \label{eq:V}\end{equation}
   In view of \e{eq:VV}   we have 
   \begin{align}
 (A\mathrm {L} f) (\xi)= &v(\xi) \sum_{m =0}^n q_m  D_{\xi}^m \big(v(\xi)^{1 -1/n} f( x(\xi))\big) 
 \nonumber\\
 =  : & v(\xi)^{-1/n} \sum_{m =0}^n i^{-m} b_{m} (x(\xi)) f^{(m)} ( x(\xi)) 
 \label{eq:V3}\end{align}
 where the functions $b_{m} (x)$ are determined by the comparison of the coefficients at $f^{(m)} ( x(\xi)) $ in the left- and right-hand sides of \e{eq:V3}. It follows that the operator  \e{eq:V} equals
     \begin{equation}
B=      \sum_{m=0}^{n} b_m (x) D^m_{x}, \q  \q D_{x}   = - id/d x , 
 \label{eq:VBb}\end{equation}
 with the coefficients $b_m (x)$ defined by equation \e{eq:V3}.
 
% Thus $B$ is also a differential operator of order $n$.
    Our goal is to find   $b_{0}(x),\ldots, b_{n}(x)$. To that end,   we need to calculate the derivatives in  the left-hand side of \e{eq:V3}. 
  For shortness, we put
  \begin{equation}
\phi (\xi) = x'(\xi) =v(\xi)^{-2/n}, \q \psi (\xi) =  v(\xi)^{1- 1/n}.
 \label{eq:V7x}\end{equation}
  By the Leibniz formula, we have
    \begin{equation}
i^{m} D_{\xi}^m \big(\psi (\xi)f ( x(\xi))\big) = \sum_{j=0}^m i^{j}\tbinom{m}{j}   \psi^{(m-j)} (\xi)    D_{\xi}^j     f( x(\xi))  .
 \label{eq:V7}\end{equation}
 The terms with $j=0$ contribute to $b_0 (x)$ which yields the expression
 \begin{equation}
 b_0 (x)=v(\xi)^{1 +1/n} \sum_{m=0}^n i ^{-m}q_m      \big(v(\xi)^{1 -1/n} \big)^{(m)}  
 \label{eq:V5}\end{equation}
 where as always the variables $x$ and $\xi$ are linked by formula  \e{eq:V2}.

To find expressions for other coefficients $b_m (x)$, we have to
  calculate the derivatives $   D_{\xi}^j     f( x(\xi))$.

% For derivatives of the function $u( x(\xi) )$, we use the following assertion.
 
    \begin{lemma}\label{Le}
    For all $j\geq 1$, we have the formula
      \begin{equation}
 i^j D_{\xi}^j f( x(\xi)) = \sum_{l=1}^j \tau_{j,l}   (\xi)   f^{(l)}( x(\xi))
 \label{eq:V8z}\end{equation}
 where
      \begin{equation}
  \tau_{j,l}   (\xi)= \sum  c_{\kappa_{1},\ldots, \kappa_{l}} \phi^{(\kappa_1)} (\xi)\cdots  \phi^{(\kappa_l)} (\xi)    ,
 \label{eq:V8}\end{equation}  
 $ c_{\kappa_{1},\ldots, \kappa_{l}}$ are some numerical coefficients and the sum is taken over all  $ \kappa_1 \geq 0, \ldots ,  \kappa_l  \geq 0$ such that
  \[
 \kappa_1 +\cdots + \kappa_l =j-l.
\]
%\label{eq:V9}\end{equation}
In particular,
 \begin{equation}
  \tau_{j,j}   (\xi)= \phi (\xi)^j
 \label{eq:V4x}\end{equation}
 and
  \begin{equation}
  \tau_{j,j-1} (\xi) = \phi'(\xi) \varphi(\xi)^{j-2}+ (\phi(\xi)^2)' \phi(\xi)^{j-3}+\cdots+ (\phi(\xi)^{j-1})' , \q j\geq 2.
 \label{eq:bm4}\end{equation}
 \end{lemma}

   \begin{pf}
   If $j=1$, then $ i D_{\xi} f ( x(\xi)) =  \phi(\xi)     f'( x(\xi))$ which is consistent with \e{eq:V8z} and \e{eq:V4x}.
   So we have to justify the passage from $j$ to $j+1$.
  Differentiating expression \e{eq:V8z}, we see that
        \[
 i^{j+1} D_{\xi}^{j+1} f( x(\xi)) = \sum_{l=1}^j \Big(  \tau_{j,l}'   (\xi)   f^{(l)}( x(\xi)) +  \tau_{j,l}   (\xi)  \phi(\xi)   f^{(l+1)}( x(\xi)) \Big)
 \]
% \label{eq:V8zz}\end{equation}
 and hence  formula \e{eq:V8z} for $ D_{\xi}^{j+1}f( x(\xi)) $  holds with
  \begin{multline}
   \tau_{j+1,1}  (\xi) =   \tau_{j,1}'   (\xi), \q
 \tau_{j+1,l}  (\xi) =   \tau_{j,l}'   (\xi)     +  \tau_{j,l-1}   (\xi)  \phi(\xi)   ,\; 2\leq l\leq j, 
  \\
  \tau_{j+1,j+1}  (\xi) =   \tau_{j,j}   (\xi)\phi(\xi) .
 \label{eq:V8zz1}\end{multline}
 Substituting here expressions \e{eq:V8} for $ \tau_{j,l}  (\xi)$,
 we get expressions \e{eq:V8} for $ \tau_{j+1,l}  (\xi)$.
 
 Formula \e{eq:V4x} follows from the last relation \e{eq:V8zz1}. Finally, relation \e{eq:V8zz1} for $l=j$ shows that
 \[
  \tau_{j+1,j}  (\xi) =    (  \phi(\xi)^j )'      +  \tau_{j,j-1}   (\xi)  \phi(\xi).
 \] 
 Using now formula \e{eq:bm4} for $ \tau_{j,j-1}  (\xi)$ we obtain   expression  \e{eq:bm4} for $  \tau_{j+1,j}  (\xi)$.
       \end{pf}

 Substituting expressions \e{eq:V7} and \e{eq:V8z} into equality \e{eq:V3} and putting together the coefficients at
$ f^{(m)} ( x(\xi)) $, we obtain the following result.

   \begin{theorem}\label{LM}
Let $A$ be  differential operator \e{eq:LOGz} where $Q(X)$ is polynomial \e{eq:LO4} with $q_{n} =1$. Under Assumption~\ref{qq} define 
the unitary operator $\mathrm{L} : L^2 (\Bbb R) \to L^2 (\Bbb R) $ by formulas \e{eq:V2}, \e{eq:V1}, and let $\phi(\xi)$,
$\psi(\xi)$  be functions  \e{eq:V7x}.    Then $B=\mathrm{L}^* A \mathrm{L} $ is  differential operator \e{eq:VBb} where   the coefficients $b_m (x)$, $m=0,1,\ldots, n$,  are sums  of terms  $($times some numerical coefficients$)$
     \begin{equation}
  v(\xi)^{1+1/n}  \phi^{(\kappa_1)} (\xi)\cdots  \phi^{(\kappa_m)} (\xi)  \psi^{(p)} (\xi), \q \xi=\xi(x) ;
 \label{eq:U}\end{equation}
 here  $ p \geq 0$,  $ \kappa_1 \geq 0, \ldots ,  \kappa_m  \geq 0$ and 
  \[
 \kappa_1 +\cdots + \kappa_m +p \leq n- m.
\]
%\label{eq:U1}\end{equation}
In particular, the coefficient $b_{n}=1$ does not depend on $x$,
 \begin{equation}
b_{n-1}(x)= q_{n-1} v(\xi (x))^{2/n}
 \label{eq:bm}\end{equation}
 and $b_{0}(x)$ is given by formula \e{eq:V5}.
 \end{theorem}
 
  Suppose now that the   polynomial $Q(X)$ is real. Then the operators $A$ and hence $B$ are symmetric. Moreover, 
  if $v(-\xi)= v(\xi)$, then $A$ commutes with the involution $\sf C$ defined by \e{eq:Complx}. Since $\mathrm{L}\sf C =\sf C \mathrm{L}$, we have the identity $B {\sf C }= {\sf C} B$ so that the coefficients of the operator $B$  satisfy the condition
    \e{eq:Compl}.
  
 \medskip
 
 {\bf 4.2.}
 Under some mild additional assumptions on the function $v(\xi)$  the coefficients  $b_m (x)$, $m=0,1,\ldots, n-1$,  of the operator $B$  constructed in Theorem~\ref{LM}  
  decay as $|x|\to \infty$.

% To show that the functions $b_m (x)$, $m=0, 1,\ldots, n- 1$,  decay as $|x|\to \infty$, we accept a  mild additional assumption.
 
  \begin{assumption}\label{qq2}
  For some $\gamma >0$ and $\d \geq 0$,
 the function  $v  $ and its derivatives satisfy estimates
   \begin{equation}
v(\xi)^{2}\leq C \big|\int_0^\xi v(\eta)^{-2/n} d\eta \big|^{- n\gamma} 
 \label{eq:HV}\end{equation}
 and    
    \begin{equation}
 | v^{(p)}(\xi)| \leq C_p v(\xi)^{1+ \d p}, \q p=1,2, \ldots.
 \label{eq:HV1}\end{equation} 
  \end{assumption}
  
   \begin{remark}\label{qr}
 Actually, it suffices to require in Assumption~\ref{qq}  that $v\in C^N (\Bbb R)$ for some sufficiently large  but finite $N $ and to impose  condition \e{eq:HV1} for $  p=1, \ldots, N$.
    \end{remark}

  Note that for function  \e{eq:LOGZ1}, condition \e{eq:HV} is satisfied with $\gamma = 1$ and estimate \e{eq:HV1}  is true for all $p$ with $\d= 0$.
   
 Putting together \e{eq:V2} and \e{eq:HV}, we see that
 \begin{equation}
v(\xi(x))^{2}\leq C (1+|x|)^{- n  \gamma }.
 \label{eq:HVv}\end{equation}
  It follows from    \e{eq:HV1}   that functions \e{eq:V7x} satisfy the estimates
     \[
    |  \phi^{(\kappa )} (\xi)| \leq C v(\xi)^{-2/n}, \q |\psi^{(p)} (\xi) | \leq C v(\xi)^{1-1/n}
    \]
    for all $\kappa$ and $p$. Therefore   expression \e{eq:U}  is bounded by  $C v(\xi)^{2-2m /n}$ and hence, by Theorem~\ref{LM},
   \[
| b_m (x)| \leq C v(\xi (x))^{2-2m /n} \leq C' (1+|x|)^{- (n -m)\gamma }, \q m =0,1,\ldots, n-1.
\]
% \label{eq:V6a}\end{equation}

 Finally, we obtain  estimates on derivatives of the functions $b_m (x)$. Note that  
    \begin{equation}
|d^p  b_m    (x)/ dx^p| \leq C_{p} \sum_{q=1}^p |d^q  b_m   (x (\xi))/ d\xi^q| \sum_{\varkappa_{j}\geq 1, \varkappa_{1 }+\cdots + \varkappa_q =p} |\xi^{(\varkappa_{1})}(x)| \cdots |\xi^{(\varkappa_q)}(x)|
 \label{eq:V6a1}\end{equation}
for all $p$. Recall that $b_m    (x (\xi))$ consists of terms  \e{eq:U}.  Differentiating them and using estimate \e{eq:HV1} we see that for all $q$,
 \[
 | d^q b_m (x (\xi))/ d\xi^q| \leq C v(\xi)^{2-2m /n + \d q} .
 \]
   Differentiating relation 
 $ \xi' (x)=  v(\xi (x))^{2/n} $
and using again estimate  \e{eq:HV1}  we find that
 \[
| \xi ^{ (\kappa)} (x)| \leq C_{\kappa}  v(\xi (x))^{2 \kappa/n + \d (\kappa-1)} .
 \]
 According to \e{eq:V6a1} this implies the bound
\[
|d^p  b_m    (x)/ dx^p| \leq C_{p} v(\xi (x))^{2-2m /n + p (\d +2/n)}. 
\]
Taking now into account condition  \e{eq:HVv}, we obtain    the following result.
   
   \begin{theorem}\label{LMb}
   In addition to the conditions of Theorem~\ref{LM}, 
    let Assumptions~\ref{qq2}    hold.
    Then $B=\mathrm{L}^* A \mathrm{L}$ is given by formula \e{eq:VBb}  where $b_{m}=1$ and the coefficients $b_{m} \in C^\infty (\Bbb R)$ 
     obey the estimates
 \begin{equation}
| b_m ^{(p)} (x)| \leq C (1+|x|)^{- (n   -m)\gamma -p \gamma (1+ \d n/2)}, \q m =0,1,\ldots, n-1, \q p=0,1,\ldots.
 \label{eq:V6ad}\end{equation}
 \end{theorem}
 
  \begin{remark}\label{GDO}
 Of course the construction above works for general differential operators
 \[
 A=\sum_{m=0}^n a_{m} (\xi) D_{\xi}^m
 \]
 under appropriate assumptions on the coefficients $a_{m} (\xi )$, $m=0,1,\ldots, n$.
   The corresponding operator $\mathrm{L}$ (cf. \e{eq:V1}, \e{eq:V2}) is defined by the formula
   \[
   (\mathrm{L} f)(\xi )= a_n (\xi)^{-1/(2n)} f\big(\int_{0}^\xi a_n (\eta)^{-1/n} d\eta\big).
   \]
 \end{remark}
 
 Let us now apply the results of Sections~2 and 3 to the operator $B$ and then reformulate them in terms of the operator $A= \mathrm{L} B \mathrm{L}^*$. As usual, we suppose that the coefficients  $q_{0},\ldots, q_{n}$ in \e{eq:VV} are real numbers. Then  the operators $A$ and $B$ are symmetric on $C_{0}^\infty(\Bbb R)$.  If
the coefficients $b_{m}(x)$, $m=0,1,\ldots, n-1$,   are bounded, then they are essentially self-adjoint. The operator $B$ is self-adjoint on $  {\sf H}^n ({\Bbb R}) $, and  the operator $A$ is self-adjoint on $\mathrm{L} {\sf H}^n ({\Bbb R}) $.  Moreover, if $b_{m}(x)\to 0$  as $|x|\to\infty$, then the essential spectra of the operators $B$ and $A$ coincide with $\Bbb R$ if $n$ is odd, and they coincide with $[0,\infty)$ if  $n\geq 2$ is even.

    If Assumptions~\ref{qq2} holds with $\gamma>1$, then according to \e{eq:V6ad} all the coefficients $b_{m} (x)$,  $m=0,1, \ldots, n- 1$,  of the operator $B$ are short-range so that   Theorem~\ref{SpTh} works  in this case.
  In our applications to Hankel operators,  only the weaker condition $\gamma>1/2$ is satisfied. Then     the coefficients $b_{m} (x)$ for $m=0,1, \ldots, n- 2$  are short-range and the coefficient $b_{n-1} (x)$ is given by formula \e{eq:bm} so that  Proposition~\ref{BL} can be used.   For an arbitrary $\gamma>0$, all the coefficients $b_{m} (x)$, $m=0,\ldots, n-1$, may be long-range, but under the assumption
 \begin{equation}
\gamma (1+ \d n/2)\geq 1
 \label{eq:V6bd}\end{equation}
  Theorem~\ref{SpThlong} applies.
    This leads to the following result.

 \begin{theorem}\label{SpThD}
 Let Assumptions~\ref{qq} and \ref{qq2}    hold. Suppose that either $\gamma>1/2$ $($and $ \d\geq 0$ is arbitrary$)$ or
condition \e{eq:V6bd} is satisfied.
Define the self-adjoint   operator $A $   by differential  expression \e{eq:LOGz} where the polynomial $Q(X)$ is  real and $q_{n}=1$.  
Then:
\begin{enumerate}[\rm(i)]
\item
The spectrum of the operator $A$ is absolutely continuous except eigenvalues that  may accumulate to zero and  infinity only.   
 \item
 Define the function $ x (\xi)$  by formula \e{eq:V2}.  
  The limiting absorption principle holds, that is, for any $\sigma>1/2$, the operator-valued  function 
  \begin{equation}
  \la x (\xi)\ra^{-\sigma}  (A -z)^{-1} \la x(\xi)\ra^{-\sigma}
 \label{eq:LAP1}\end{equation}
 is continuous up to the real axis,  except the eigenvalues of the operator $A$ and the point zero.
   \item
  The 
 absolutely continuous  spectrum of the operator $A$ covers $\Bbb R$ and is simple for $n$ odd. It coincides with $[0,\infty)$ and has multiplicity $2$  for $n$ even. 
  \item
 If $n$ is odd, then the multiplicities of eigenvalues  of the operator $A$  are bounded by $(n-1)/2$. If $n$ is even, then the multiplicities of positive eigenvalues are bounded by $n/2-1$, and the multiplicities of negative eigenvalues are bounded by $n/2$. 
   \end{enumerate} 
      \end{theorem}

We emphasize that since the function $x(\xi)/\xi \to \infty$ as $|\xi|\to \infty$,   the limiting absorption principle
for the   operator $A$ with degenerate coefficients  has a weaker form than in the regular case discussed in Sections~2 and 3. In particular, for function  \e{eq:LOGZ1} the resolvent $(A -z)^{-1} $ has to be sandwiched by the exponentially decaying weight $  \la x (\xi)\ra^{-\sigma}  $.

 \begin{remark}\label{GDOx}
 Condition \e{eq:V6bd} is far from being optimal. It is required only to use conveniently the results of Section~3 on the long-range case (Theorem~\ref{SpThlong}).
  \end{remark}

% In the  case $\gamma>1/2$ estimates \e{eq:V6ad}  mean that the coefficients $b_{m}(x)$, $m=0,1,\ldots, n-2$,
% satisfy condition \e{eq:SR} with $\rho=2\gamma > 1$ and $b_{n-1}(x)$ satisfy condition \e{eq:LR} with $\rho=\gamma > 1/2$. Therefore making the unitary transformation \e{eq:Ga} we can reduce the operator $B$ to the operator $\wt{B}$ with short-range coefficients.   

Let us discuss the case $\gamma>1/2$ in more details. Following Proposition~\ref{BL}, we define    eigenfunctions $\psi_{\pm}(x,k)$ of the operator $B$ by formula \e{eq:Ga3} where ${\wt\psi}_{\pm}(x,k)$ are the eigenfunctions of 
  the differential operator $\wt B ={\cal J}^* B {\cal J}$
with short-range coefficients. It follows from formulas \e{eq:Ga1} and \e{eq:bm} that 
 \begin{equation}
\beta (x)= -n^{-1}q_{n-1} \int_{0}^x v(\xi(y))^{2/n}  dy= - n^{-1}q_{n-1} \xi (x)
 \label{eq:gauge}\end{equation}
 where we have used that $dy= v(\eta)^{-2/n}d\eta$. Thus relations  \e{eq:Ga} and  \e{eq:Ga3} read now as
  \begin{equation}
 ( {\cal J}  f)(x) = e^{- i q_{n-1} \xi (x) /n} f(x)
\label{eq:Gab}\end{equation}
and
  \begin{equation}
  \psi_{\pm}(x,k)= e^{-i  q_{n-1} \xi (x)/n} {\wt\psi}_{\pm}(x,k).
 \label{eq:gauge1}\end{equation}

According to \e{eq:V} we define eigenfunctions of the operator $A$ by the formula $\varphi_{\pm}( k) = \mathrm{L} \psi_\pm(k)$. Then representation \e{eq:est2B}  implies that the spectral family $E_{A} (\lambda)$ of the operator $A$
 satisfies the relation
 \[
\frac{d(E_{A}(\lambda)g,g)}{d\lambda}=\sum_{k^n=\lambda}\big|\int_{-\infty}^\infty   \ov{\varphi_{\pm} (\xi,k)} g(\xi) d\xi\big|^2, \q g\in C_{0}^\infty ({\Bbb R} ),  \q \lambda=k^n\in \Omega.
\]
%\label{eq:est2A}\end{equation}
 It follows from
\e{eq:V2}, \e{eq:V1} and \e{eq:gauge1} that 
  \begin{equation}
\varphi_{\pm}(\xi,  k) = e^{- iq_{n-1}\xi /n} v(\xi)^{-1/n} \wt\psi_{\pm} \big(\int_{0}^\xi v(\eta)^{-2/n} d\eta, k \big).
 \label{eq:eig}\end{equation}

%Here the sum consists of the one term if $n$ is odd; if $n$ is even it consists of two terms for $\lambda> 0$, and it is empty if $\lambda<0$.

 Using the results of Section~2,      we can write down asymptotics 
 of the functions $\varphi (\xi,  k) = \varphi_{-}(\xi,  k)$ as $\xi\to \pm\infty$. For example,
for odd $n$, it follows from \e{eq:WFSodd2} and \e{eq:eig} that
%\begin{equation}
\[
\left\{\begin{array}{lcl} 
\varphi (\xi,  k) = e^{- iq_{n-1}\xi /n} v(\xi)^{-1/n}  \exp\big( ik\int_{0}^\xi v(\eta)^{-2/n} d\eta\big) (1+o(1)) ,\quad \xi\ri -\infty,
   \\
\varphi(\xi,  k) =   s(\lambda) e^{- iq_{n-1}\xi /n} v(\xi)^{-1/n}  \exp\big( ik\int_{0}^\xi v(\eta)^{-2/n} d\eta\big) (1+o(1)) ,\quad \xi\ri \infty,
\end{array}\right.
\]
%\label{eq:eig1}\end{equation}
where $s(\lambda)$ is the scattering matrix for the pair $B_{0} , \wt B$ (or which is the same, for the triple $B_{0} , B$, $\cal J$). Quite similarly, for even $n$ the asymptotics of 
$\varphi (\xi,  k) $ follows from 
\e{eq:WFS1}, \e{eq:WFS2}  and \e{eq:eig}. We recall that the asymptotic coefficients $s_{j \ell} (\lambda)$ in these formulas are elements of the scattering matrix $S(\lambda)$ for the pair  $B_{0}$, $\wt B$. In particular, if $v(-\xi)= v(\xi)$, then the operator $A$ commutes with the involution $\cal C$ defined by \e{eq:Complx}. In this case we also have ${\cal C} \mathrm{L}=\mathrm{L}{\cal C}$ so that ${\cal C} B= B {\cal C}$. Thus $S(\lambda)$  satisfies identity \e{eq:Complz}, and hence $s_{12} (\lambda)= s_{21} (\lambda)$.

 The eigenfunctions \e{eq:eig} grow and rapidly oscillate as $|\xi|\to\infty$. For example, for function
   \e{eq:LOGZ1} both the amplitude and the phase of $\varphi (\xi,  k) $ grow
    exponentially. This is of course consistent with the formulation of the limiting absorption principle in part (ii) of Theorem~\ref{SpThD}.

\medskip

{\bf 4.3.}
 Let us discuss conditions  \e{eq:HV} and \e{eq:HV1}.

 \begin{example}\label{exp1}
Suppose that
  \[
c_1 (1+  |\xi|)^{-\alpha} \leq v(\xi)\leq c_{2}  (1+  |\xi|)^{-\alpha}, \q   \alpha>0, \q c_{1}, c_{2}  >0.
 \]
 % \label{eq:exp1}\end{equation}
 Then 
  \begin{equation}
\Big| \int_0^\xi v(\eta)^{-2/n} d\eta \Big| \leq c_1^{-2/n} \Big|\int_0^\xi(1+ | \eta |)^{2\alpha/n}  d\eta \Big| \leq  C_{1}  (1+ |  \xi |)^{1+ 2\alpha/n} 
 \leq  C_2  v(\xi )^{- (n+ 2\alpha)/ (\alpha n)} 
 \label{eq:exp1a}\end{equation}
so that condition  \e{eq:HV}  is satisfied with $\gamma=2\alpha (2\alpha+n)^{-1}$.

  If, moreover, 
   \[
|  v^{(p)}(\xi)| \leq C_{p}  (1+  |\xi|)^{-\alpha-p}
 \]
 for all $p=0,1,\ldots$, then condition \e{eq:HV1} is satisfied with $\d=1/\alpha$. Therefore we have the equality in  \e{eq:V6bd} and estimate \e{eq:V6ad}  reads as 
   \[
| b_m^{(p)} (x)| \leq C (1+|x|)^{- (n -m)\gamma -p  }, \q m =0,1,\ldots, n-1, \q \gamma = 2\alpha (2\alpha +n)^{-1}.
 \]
 Thus the coefficients $b_m (x)$ decay as $|x|\to \infty$, but at least some of them are long-range.
 
 It follows from \e{eq:exp1a} that function \e{eq:V2}  satisfies the estimate 
$
 \la x( \xi) \ra \leq C  \la   \xi \ra^{1+ 2\alpha/n}.
$
 Therefore in the formulation of part (ii) in Theorem~\ref{SpThD}, the operator-valued  function \e{eq:LAP1} can be replaced by
 $
  \la  \xi\ra^{-\sigma}  (A -z)^{-1} \la  \xi\ra^{-\sigma}
$
 where $\sigma>1/2 + \alpha/n$.
 \end{example}

  \begin{example}\label{exp}
Suppose that
  \[
c_1 e^{ - \beta |\xi|^\alpha}\leq v(\xi)\leq c_{2} e^{ - \beta |\xi|^\alpha}, \q  \alpha>0, \;  \beta >0.
\]
% \label{eq:exp}\end{equation}
 Then 
  \begin{equation}
 \Big|\int_0^\xi v(\eta)^{-2/n} d\eta  \Big| \leq c_1^{-2/n}  \int_{0}^{|\xi|} e^{ 2 \beta  \eta ^\alpha /n} d\eta\leq  C_{1} e^{ 2 \beta | \xi|^\alpha /n} (1+ |\xi|)^{1-\alpha}.
 \label{eq:expa}\end{equation}
so that condition  \e{eq:HV}  is satisfied with $\gamma=1$ if $\alpha\geq 1$ and with arbitrary  $\gamma<1$ if $\alpha\in(0, 1)$.   If, moreover, 
   \[
 | v^{(p)}(\xi)| \leq C_{p}e^{ - \beta |\xi|^\alpha}
 \]
 for all $p=0,1,\ldots$, then condition \e{eq:HV1} is satisfied with $\d=0$.
 
 In particular, for function  \e{eq:LOGZ1} $\alpha=1$, $\beta=\pi/2$, and hence according to \e{eq:V6ad}    the coefficients of the corresponding operator $B$ satisfy for all $p$ the estimates 
 \[
| b_m^{(p)} (x)| \leq C (1+|x|)^{-n+m -p}, \q m=0,1,\ldots, n-1 .
\]
% \label{eq:expV}\end{equation}
 In this case the coefficients $b_{0}(x),\ldots, b_{n-2} (x)$ are short-range and  $b_{n-1}^{(p)} (x) =O(|x|^{- 1 -p})$. 
 
 It follows from \e{eq:expa} that in the formulation of the limiting absorption principle in Theorem~\ref{SpThD} the operator-valued  function \e{eq:LAP1} can be replaced by
 $
  \la   e^{  |\xi|^\alpha}\ra^{-\sigma}  (A -z)^{-1}  \la  e^{  |\xi|^\alpha}\ra^{-\sigma}   
$
 where $\sigma> \beta/n$.
 \end{example}

 % Maybe to distinguish the notation for the variable $x$ and the function $\xi$?

   %************************************************************
\section{Hankel operators. Spectral results} 
%********************************************************* 

Here we come back to the generalized Carleman operators $H$. Our goal in this section is to prove Theorem~\ref{SpThHa}.

\medskip
 
  {\bf 5.1.}
So we consider  Hankel operators  
 \begin{equation}
(H u)(t) = \int_{0}^\infty h(t+s) u(s)ds 
\label{eq:H1bis}\end{equation}
  with kernels 
   \begin{equation}
h(t)= t^{-1} \sum_{m=0}^n p_m \ln^m t ,\q p_{m} =\bar{p}_{m}, \q p_{n}=1, \q n\geq 1.
\label{eq:LOGbis}\end{equation} 
%The coefficients $p_m$, $m=0,1,\ldots, n$, are supposed to be real and $p_{n}=1$.
Put
 \begin{equation}
q_m=    \sum_{j =m}^n  \tbinom{j}{m} \gamma^{(j-m)} (0) p_j ,\q m=0,\ldots, n, \q \gamma ( z)=\Gamma (1-z)^{-1},
\label{eq:LOG5a}\end{equation} 
where $\Gamma(\cdot)$ is the gamma function.  
For example,
$q_{n}= p_n $ and
 \[
 q_{n-1}= p_{n-1}  + \Gamma' (1) n  p_n
  \]
  %  \label{eq:LOG61}\end{equation} 
    for all $n$ (recall that  $  -\Gamma' (1)$ is the Euler constant).  Of course formulas \e{eq:LOG5a} allow one to recover the coefficients $p_{n}, p_{n-1},\ldots, p_{0}$ given the coefficients $q_{n}, q_{n-1},\ldots, q_0$.

 Let $A$ be the  differential operator 
   \begin{equation}
A= v( \xi)
\Big(\sum_{ m=0}^n q_m D_{\xi}^m \Big) v( \xi),\q D_{\xi}   = -i d /d \xi,
\q v(\xi)=\frac{\sqrt{\pi}} {\sqrt{\cosh (\pi \xi)}},
\label{eq:LOGA}\end{equation}
 in the space $L^2 ({\Bbb R}) $.
It is self-adjoint on the domain ${\cal D}(A)=\mathrm{L} {\sf H}^{n}({\Bbb R})$ where the unitary operator $\mathrm{L}$ is defined by formula  \e{eq:V1} with $x(\xi)$ given by
    \begin{equation}
x(\xi)=\pi^{-1/n}\int_{0}^\xi \big(\cosh( \pi s) \big)^{1/n} d s.
 \label{eq:han}\end{equation}

The  operators $H$ and $A$ turn out to be unitary equivalent. Let $M: L^2 ({\Bbb R}_{+}) \to L^2 ({\Bbb R} )$ be the Mellin transform defined by the formula
 \[
(M u) (\xi)=  (2\pi)^{-1/2} \int_0 ^\infty u(t) t^{-1/2 -i \xi  } dt .
\]
It is a unitary mapping. Set
  \begin{equation}
(F u) (\xi)= e^{ -i\eta(\xi)} (M u)(-\xi)\q{\rm where}\q e^{- i\eta(\xi)} =\frac{\Gamma(1/2- i \xi)}{|\Gamma(1/2 - i \xi)|}.
\label{eq:MAID1}\end{equation}
We proceed from the following result. Here we note only that its proof relies on the identity
\[
h(t)=\int_{0}^\infty e^{-\lambda t} \sum_{m=0}^n q_{m} \ln^m \lambda d\lambda.
\]

 \begin{theorem}\label{1z} \cite[Theorem~3.2]{Yf1}
 Let $H$  be the Hankel operator \e{eq:H1bis} with  kernel  given  by   \e{eq:LOGbis}. Let  the coefficients $q_{m}$ be defined by formulas \e{eq:LOG5a}, and let $A$ be   differential operator \e{eq:LOGA}. Then for all functions $u_{j} $, $j=1,2$, such that
 $M u_{j}\in C_{0}^\infty ({\Bbb R})$, 
 the identity 
 \begin{equation}
(H u_{1}, u_{2})=  (A F u_{1}, F u_{2})
 \label{eq:MM1}\end{equation} 
holds.
 \end{theorem}

  We emphasize that Theorem~\ref{1z} does not require that the coefficients of $P(X)$ be real. Note that for $n=0$, Theorem~\ref{1z} reduces to the diagonalization of the Carleman operator. If $n\geq 1$, then the operators $H$ are unbounded, but for real coefficients $p_{m}$ they are self-adjoint on the set ${\cal D}(H)=F^{*}  {\cal D}(A)$. In this case 
  $A {\sf C } ={\sf C }A  $ where $\sf C$ is  the involution  given by \e{eq:Complx}.  
  
  For  the function $v(\xi)$ defined in  \e{eq:LOGA}, conditions \e{eq:HV} and \e{eq:HV1} are satisfied  with $\gamma=1$ and $\d=0$,  respectively. Thus Theorem~\ref{SpThD}  implies the spectral results (i), (iii) and (iv) of Theorem~\ref{SpThHa}. 
  
 \medskip

%  Our goal is to find asymptotics of these functions as $t\to \infty$ and  as $t\to 0$, but we do not dwell upon the precise justifications of all formulas.

  {\bf 5.2.}
    The limiting absorption principle   requires a special discussion. By virtue of part (ii) of Theorem~\ref{SpThD}  and the identity \e{eq:MM1}        the operator-valued  function 
 \[
  \la x (\xi)\ra^{-\sigma} M (H -z)^{-1} M^* \la x(\xi)\ra^{-\sigma}, \q \sigma>1/2,
 \]
 % \label{eq:LAPh}\end{equation}
   is continuous up to the real axis in the complex plane  cut along $\Bbb R$ for $n$ odd and along $[0,\infty)$ for $n$ even,  except the eigenvalues of the operator $H$ and the point zero. However this does {\it not}  imply the same conclusions about  the operator-valued  function \e{eq:LAPH}. 
   
   To prove part (ii) of Theorem~\ref{SpThHa},  we first diagonalize  the   Hankel operator $H$. Comparing   relations \e{eq:F2}, \e{eq:V}  and \e{eq:MM1}, we   see that
    \begin{equation}
\Theta_{\pm} H = \Lambda\Theta_{\pm} \q {\rm  where} \q \Theta_{\pm} =\Psi_{\pm} \mathrm{L}^* F,
 \label{eq:Diag}\end{equation}
the operator  $\mathrm{L}$ is defined by relations \e{eq:V2} and, as before, $ \Lambda$ is the operator of multiplication by $k^n$ in the space $L^2 ({\Bbb R})$.
The operators $\Psi_{\pm} $     are given by relation \e{eq:WFsr1} in terms of the eigenfunctions ${   \psi}_{\pm}(x,k)$  of the continuous spectrum  of the operator $B=\mathrm{L}^* A \mathrm{L}$.
 Recall that,  with the operator $\cal J$ given by  \e{eq:Gab},   the differential operator $\wt B={\cal J}^*  B {\cal J}$ has short-range coefficients. Therefore its eigenfunctions ${ \wt \psi}_{\pm} (x,k)$ are correctly defined by formula \e{eq:WFt} and 
 the eigenfunctions ${   \psi}_{\pm}(x,k)$ of the operator $B$ are given by \e{eq:gauge1}.

It follows from   \e{eq:Diag} that eigenfunctions $\theta_{\pm} (t,k)$ of the continuous spectrum of the operator $H$ can formally be defined by the equation 
$
\theta_{\pm} (k)= F^* \mathrm{L} \psi_{\pm} (k)
$
(note  that $\varphi_{\pm} (k)=   \mathrm{L} \psi_{\pm} (k)$ are the   eigenfunctions of the operator $A$)
or, in a more detailed notation,  
 \begin{align}
\theta_{\pm} (t,k)=&  (2\pi)^{-1/2} t^{-1/2 } \int_{-\infty}^\infty e^{  -i\xi \ln t}
e^{i \eta(\xi)} \varphi_{\pm} (\xi,k) d\xi
\nonumber\\
=&   (2\pi)^{-1/2} t^{-1/2 } \int_{-\infty}^\infty e^{  -i\xi (x)\ln t}
e^{i \eta(\xi (x))} \xi'(x)^{1/2}   \psi_{\pm}  (x,k) d x
 \label{eq:han1X}\end{align}
 where $\xi (x)$ is the function inverse to \e{eq:han}.
As we will see in the next section, these  integrals  converge (but not absolutely). 
Of course, we always suppose that $\lambda=k^{n} $ belongs to the set $\Omega$ defined by relations \e{eq:omega}. 
  By their construction, the functions $\theta_{\pm} (t,k)$   satisfy the equation $H \theta_{\pm}(k)=k^{n} \theta_{\pm}(k)$. In terms of the functions $\theta_{\pm} (t,k)$, we have (cf. \e{eq:WFsr1})
 \begin{equation}
({\Theta}_{\pm} u)(k)= (2\pi)^{-1/2}  
 \int_{0}^\infty \overline{\theta_{\pm}(t, k)}u(t) dt. 
 \label{eq:HAN}\end{equation}

In the next section, we will establish the following result. As usual, we set $\theta=\theta_{-}$.

   \begin{theorem}\label{EIG} 
   The functions $\theta (t,k)$ defined by formula \e{eq:han1X} satisfy the estimates 
    \begin{equation}
    | \theta (t , k)|\leq C t^{-1/2}
     \label{eq:est}\end{equation}
  and   
      \begin{equation} 
|\partial\theta(t,k)/ \partial k| \leq C t^{-1/2} \la\ln t\ra .
\label{eq:est1dif}\end{equation}
These estimates are uniform in $k$ for $\lambda=k^{n} $ in compact subsets of $\Omega$.
 \end{theorem}

We note that estimates \e{eq:est} and \e{eq:est1dif} on the eigenfunctions $\theta(t,k)$ of the operator $H$ are the same  as for the classical Carleman operator when $h(t)=t^{-1}$.

    \begin{corollary}\label{smoothH} 
    For any $\alpha\in (0,1)$ and $\sigma>\alpha+1/2$ the operator 
    \[
     \Theta \la\ln t\ra^{-\sigma}: L^2
 \to C^\alpha (\Omega)
   \]
   %  \label{eq:Smooth}\end{equation}
is bounded.
 \end{corollary}
 
   \begin{pf}
   It follows from definition \e{eq:HAN} and estimate \e{eq:est} that
       \begin{equation} 
   |( \Theta \la\ln t\ra^{-\sigma} u) (k)|^2 \leq (2\pi)^{-1} \big( \int_{0}^\infty |\theta (t, k)| \la\ln t\ra^{-\sigma}|u(t) |dt\big)^2
   \leq C \int_{0}^\infty t^{-1} \la\ln t\ra^{-2\sigma} dt \| u\|^2.
  \label{eq:Smooth1}\end{equation}
   Hence, for $\sigma>1/2$, the functions $( \Theta \la\ln t\ra^{-\sigma} u) (k)$ are bounded  uniformly in $k$  in compact subsets of $\Omega$.
 
 Comparing estimates \e{eq:est} and \e{eq:est1dif}, we see that
  \[ 
|\theta(t,k)-\theta(t,k')| \leq C t^{-1/2} \la\ln t\ra^\alpha |k-k'|^\alpha, \q \forall \alpha\in (0,1).
\]
%\label{eq:est1}\end{equation}
Therefore quite similarly to \e{eq:Smooth1}, we find that
\[
   |( \Theta \la\ln t\ra^{-\sigma} u) (k) -( \Theta \la\ln t\ra^{-\sigma} u) (k')|  \leq C |k-k'|^\alpha  \big(  \int_{0}^\infty t^{-1} \la\ln t\ra^{-2\sigma+ 2 \alpha} dt \big)^{1/2} \| u\|.
   \]
The last integral converges if $\sigma>\alpha+1/2$.
   \end{pf}
   
In the standard terminology of scattering theory (see, e.g., Definition~4.5 in Chapter 4 of the book \cite{YaMSC}) Corollary~\ref{smoothH} means that the operator $ \la\ln t\ra^{-\sigma}$ is strongly $H$-smooth with exponent $\alpha$ on $\Omega$. In particular,  the operator $ \la\ln t\ra^{-\sigma}$ for $\sigma>1/2$ satisfies Definition~\ref{Kato}.

It follows from formula \e{eq:Diag} that the spectral family  $E_{H} (\lambda)$ of the operator  $H$ admits    
   the representation
 \[
\frac{d(E_{H} (\lambda) u,u )}{d\lambda}= n^{-1} k^{-n+1}\sum_{k^n=\lambda} | (\Theta u)(k)|^2, \q \lambda\in \Omega.
\]
%\label{eq:est2}\end{equation}
Here the sum consists of the one term if $n$ is odd; if $n$ is even it consists of two terms for $\lambda> 0$, and it is empty if $\lambda<0$. Therefore Corollary~\ref{smoothH} implies that the operator-valued function
\[
d \la \ln t \ra^{-\sigma} E_{H} (\lambda)\la \ln t \ra^{-\sigma}/d\lambda, \q \sigma >1/2,
\]
is H\"older continuous with exponent $\alpha<\sigma-1/2$ in $\lambda\in \Omega$. Using now the Privalov theorem (see, e.g., Theorem~2.6 in Chapter~1 of \cite{YaMSC}) we conclude the proof of statement (ii) of Theorem~\ref{SpThHa}.

%  Theorem~\ref{EIG} allows us to  prove the limiting absorption principle for the operator $H$. 

     \medskip
 
  {\bf 5.3.}
   The limiting absorption principle for the operator $H$ allows one to easily extend Theorem~\ref{SpThHa} to perturbations of generalized Hankel operators. We state only the simplest results in this direction.  
  However we do not assume that a perturbation ${\bf V}$ is a Hankel operator. Instead, we accept a more general
  
    \begin{assumption}\label{VV}
An    operator ${\bf V}$ is self-adjoint and the operator $\la \ln t\ra^{\alpha_{1}} {\bf V} \la  \ln t\ra^{\alpha_{2}}$ is compact for some $\alpha_{1}, \alpha_{2}>1/2$.
    \end{assumption}
    
    Recall that the wave operators were introduced in subs.~2.2.
  
    \begin{theorem}\label{Perturb}
 Let $H$  be the Hankel operator \e{eq:H1bis} with  kernel  given  by   \e{eq:LOGbis}. Under Assumption~\ref{VV} put ${\bf H}=H+{\bf V}$.  Then 
\begin{enumerate}[\rm(i)]
 \item
The wave operators $W_{\pm} ({\bf H}, H)$ exist and are complete. In particular, the 
 absolutely continuous  spectrum of the operator $\bf H$ covers $\Bbb R$ and is simple for $n$ odd. It coincides with $[0,\infty)$ and has multiplicity $2$  for $n$ even. 
 \item
 If moreover $\alpha_{j}>1$ for at least one $j$, then the operator $\bf H$ has no singular continuous spectrum. The eigenvalues of $\bf H$ distinct from $0$ have finite multiplicities and can accumulate to this point $($and infinity$)$ only.
\end{enumerate}
    \end{theorem}
    
   The proof of this result relies on the tools of scattering theory.  We can use, for example, Theorem~6.4 for the proof of (i) and  Theorem~7.9, 7.10 in Chapter~4 of the book \cite{YaMSC} for the proof of (ii).
   We omit details moreover that the proof of Theorems~\ref{Perturb} is practically the same as the proof of Theorems~5.1 and 5.4 in \cite{Y2} where the  case $h(t)=t^{-1}$ was considered.
    
    If $\bf V$ is an integral operator,
    \[
    ({\bf V} u)(t)=\int_{0}^{\infty} {\bf v}(t,s) u(s) ds,
    \]
    then the operator $\la \ln t\ra^{\alpha_{1}} {\bf V} \la  \ln t\ra^{\alpha_{2}}$ belongs to the Hilbert-Schmidt class and hence it is compact if 
     \[
\int_{0}^{\infty} \int_{0}^{\infty}  | {\bf v}(t,s)|^{2} \la \ln t\ra^{2\alpha_{1}} \la \ln s\ra^{2\alpha_{2}} dt ds<\infty.
\]
% \label{eq:pert}\end{equation}
 In particular, for Hankel operators ${\bf V}$ when ${\bf v}(t,s)=  v(t+s)$, this condition is satisfied if
  \[
  \int_{0}^{\infty}  |  v(t)|^{2} \la \ln t\ra^{2 \alpha }  tdt  <\infty
 \]
 %\label{eq:pert1}\end{equation}
 where $\alpha=\alpha_{1}+\alpha_{2}$. Thus assertion (i) of Theorem~\ref{Perturb} is true if $\alpha>1 $, and assertion (ii)   is true if $\alpha>3/2 $.

    %************************************************************
\section{Eigenfunctions of Hankel operators} 
%********************************************************* 

      In this section we  find the asymptotic behavior   as $t \to\infty$ and as $t \to 0$ of the eigenfunctions  $\theta (t,k)$ of the Hankel operators $H$  defined by formulas \e{eq:H1bis} and \e{eq:LOGbis}.  In particular, we prove   Theorem~\ref{EIG} and hence conclude the proof of Theorem~\ref{SpThHa}.
     
      \medskip
 
  {\bf 6.1.}
  Recall that the eigenfunctions $\theta (t,k)=\theta_{-} (t,k)$  of the Hankel operator $H$ were defined by equality  \e{eq:han1X} where $\psi (x,k)$ are the eigenfunctions of the differential operator $  B= \mathrm{L}^* F H F^* \mathrm{L}$. Taking into account formula \e{eq:gauge1} we see that
     \begin{equation}
\theta (t,k) 
=   (2\pi)^{-1/2} t^{-1/2 } \int_{-\infty}^\infty e^{  -i\xi (x)\ln t}
e^{i  \varrho (x)} \xi'(x)^{1/2}   \wt\psi  (x,k) d x
 \label{eq:han1}
 \end{equation}
 where  
 \begin{equation}
\varrho(x)= \eta (\xi (x))- q_{n-1}  \xi (x)/n
\label{eq:han4}\end{equation}
and $\wt\psi  (x,k)$ are the eigenfunctions of
  the differential operator $\wt B={\cal J}^{*} B{\cal J}$ with short-range coefficients.  We denote by $s(\lambda)$  and $s_{j \ell}(\lambda)$, $j , \ell =1,2$, $\lambda =k^{n}\in \Omega$,   the elements of the scattering matrices $($see subsection~$2.3)$ for the operators $B_{0} =D^n$, $\wt B $.  Note that $|s(\lambda)|=1$ for odd  $n$. If $n$ is even, then the scattering matrix \e{eq:SM} with the elements $s_{j \ell}(\lambda)$ is unitary and     $s_{12}(\lambda)=s_{21} (\lambda)$ according to Proposition~\ref{BL1}.
  
To give a precise definition of the phase function $\omega(t,k)$ in \e{eq:omeg}, we need some notation. For $s>0$, let us set
   \begin{equation}
\gamma_{0} (s)= -  n \pi^{-1}  \ln ( (2\pi)^{1/n}s) +n \pi^{-1},  
 \label{eq:STP2q}\end{equation}
  \[
\gamma_{1}(s)=\pi^{-1}   \ln ( (2\pi)^{1/n} s) \Big(n \ln\big|\pi^{-1} n \ln ( (2\pi)^{1/n} s)\big|-n -q_{n-1}\Big) 
\]
% \label{eq:BFG}\end{equation}
and
 \begin{equation}
\gamma (N,k)=  N \gamma_{0} (| N /k|)+ \gamma_1 (| N/k|) +\sgn N \, ( \pi/4 + a_{1}|k|) 
 \label{eq:BFGT}\end{equation}
 where
     \begin{equation}
 a_{1} =   (2\pi)^{-1/n} \int_{0}^\infty \big( (2\cosh (\pi\xi))^{1/n}  
 -e^{\pi\xi/n} \big) d\xi - \pi^{-1} (2\pi)^{-1/n}  n .
 \label{eq:BH4}\end{equation}
  
 Now we are in a position to state the main result  
of this section.

 \begin{theorem}\label{SpThHb}
Let $\theta (t,k)$ be the eigenfunctions      of the Hankel operator \e{eq:H1bis} with  kernel  given  by   \e{eq:LOGbis}.
%Let    $   \wt \psi (k)$ be eigenfunctions of the operator   $\wt B={\cal J}^{*}U^* F H F^* U{\cal J}$, and let    $  u (k) =F^* U  {\cal J}\wt\psi (k)$ be eigenfunctions of the operator   $H$.   Let the set $\Omega$ be defined by \e{eq:omega}.
 Define the phase function $\omega (t,k)$ by the formula
 \begin{equation}
\omega (t,k)=   \gamma (\ln t,k) .
 \label{eq:BFGt}\end{equation}
  Then: 
  \begin{enumerate}[\rm(i)]
\item
For $n$   odd,   relation \e{eq:Sodd2H} holds as $t\to\infty$ if $k>0$ and as $t\to 0$ if $k<0$. If $t\to 0$ and  $k>0$ or $t\to \infty$ and $k>0$, then  estimate \e{eq:So-} is satisfied.
 \item
For $n$ even,  relations \e{eq:S1H} hold if $k>0$ and relations \e{eq:S2H} hold if $k<0$. 
 \end{enumerate}
 All these relations  are uniform in $\lambda=k^n$ in compact subsets of $\Omega$.
   \end{theorem}

Of course, asymptotic formulas \e{eq:Sodd2H}, \e{eq:So-} and \e{eq:S1H}, \e{eq:S2H} imply 
  estimate  \e{eq:est}. To obtain  estimate \e{eq:est1dif}, we formally differentiate the asymptotic formulas for $\theta(t,k)$ in $k$ and observe   that $|\partial \omega(t,k)/\partial k|\leq C (1+ | \ln t|)$. To give the precise proof, we have to differentiate the integral representation  \e{eq:han1} and take into account that $\partial \psi(x,k)/\partial k $ contains an additional factor $x$ compared to $\psi(x,k)$. Then we can repeat the same calculation as for the integral \e{eq:han1} itself. We will not dwell upon details. Thus given  Theorem~\ref{SpThHb}, we conclude the proof of  Theorem~\ref{EIG} and hence of  Theorem~\ref{SpThHa}.

\medskip 
 
{\bf 6.2}.
For the proof of   Theorem~\ref{SpThHb}, we have to
    find the asymptotic behavior  as $t\to \infty$ and as $t\to 0$ of the integral \e{eq:han1}.
    % where  $\wt \psi(x,k)$ are the eigenfunctions of the operator $\wt B$ with short-range coefficients. Below we use freely the results of Section~2 and, to simplify notation, we remove tilde    from notation of various objects related to the operator $\wt B$. 
    Recall that      the functions $x(\xi)$ and $\eta(\xi)$ were defined by formulas \e{eq:han} and \e{eq:MAID1}, respectively;  $ \xi (x)$ is the function inverse to $x(\xi)$.         Put $N=\ln t$,
\begin{equation}
\zeta (x)=   e^{i \varrho (x)}    \xi'(x)^{1/2}  
 \label{eq:BH1}\end{equation}
  and 
\begin{equation}
{\cal I} (N,k)=     \int_{-\infty}^\infty e^{  -i\xi (x) N  } \zeta (x)  \wt\psi(x,k) dx
 \label{eq:BH}\end{equation}
 where   $ \wt\psi(x,k)$ are the eigenfunctions of the operator $\wt B$.
  Then according to \e{eq:han1}   we have
 \begin{equation}
 \theta(t,k)=(2\pi)^{-1/2} t^{-1/2} {\cal I}  (\ln t,k).
 \label{eq:BHC}\end{equation}
 Thus our goal is to find the asymptotics of the integral \e{eq:BH} as $N\to\pm\infty$.

By definition \e{eq:han},  the  function $x(\xi)$  is odd, $x'(\xi)>0$ and 
  \begin{equation}
x(\xi)= a_{0}^{-1}   e^{\pi\xi/n} + a_{1} + O( e^{- \pi (2n-1)\xi/n}), \q \xi\to+\infty,
 \label{eq:BH3}\end{equation}
 where  $a_{0}= \pi  (2\pi)^{1/n}  n^{-1} $ and $ a_{1} $ is given by formula \e{eq:BH4}.
 The inverse function $\xi(x)$  is also odd, $\xi'(x)>0$, and
 it follows from \e{eq:BH3} that  
   \begin{align}
\xi(x)&= \pi^{-1}    n \big( \ln ( a_{0}x) -    a_{1}x^{-1} + O( x^{- 2})\big),  
 \label{eq:BH5}\\
\xi'(x)&= \pi^{-1}    n  x^{-1}\big( 1 +   a_{1}x^{-1} + O( x^{- 2})\big), 
 \label{eq:BH6}\\
\xi''(x)&=- \pi^{-1}    n  x^{-2}\big( 1 +     O( x^{- 1})\big) 
 \label{eq:BH7}\end{align}
 as $x\to+\infty$.
  The  last relation can be further differentiated.
 
According to the Stirling formula the function $\eta(\xi)=- \eta(-\xi)$ defined by  \e{eq:MAID1}   satisfies the asymptotic relation
 \[
 \eta(\xi)= \xi \ln|\xi| -\xi + O (\xi^{-1}), \q |\xi |\to \infty.
 \]
 Therefore \e{eq:BH5} implies the  asymptotic  formula for the function \e{eq:han4}:
 \begin{equation}
\varrho(x)= \pi^{-1}   \ln (a_{0} x)\Big(n \ln \big|\pi^{-1}  n \ln (a_{0} x) \big| -n -q_{n-1}\Big)+ O( x^{-1}), \q x\to +\infty.
\label{eq:han4q}\end{equation}
 This formula can be   differentiated.
Of course the function $\varrho(x)$ is   odd.

 It follows from \e{eq:BH6}, \e{eq:BH7} and \e{eq:han4q} that function \e{eq:BH1} satisfies the estimates
     \begin{equation} 
    \zeta(x)=O(|x|^{-1/2 }), \q  \zeta'(x)=O(|x|^{-3/2 }\ln\ln |x|) 
      \label{eq:BHNd1}\end{equation}
       as $|x|\to\infty$.

%     The stationary points of such integrals are determined by the equation
%  $ \xi'(x)\ln t=k $   for $x$.   Since    $ \xi'(x)>0$, this equation  has no solutions if $t\to\infty$ and $k<0$ or if $t\to 0$ and $k>0$. On the contrary, it has  two solutions (a positive and a negative)  if $t\to\infty$ and $k>0$ or if $t\to 0$ and $k<0$. 

     We proceed from decomposition \e{eq:LS}  of the function $ \wt\psi(x,k)  $ into oscillating $\psi_{\rm osc}  (x,k) $ and decaying  $\psi_{\rm dec} (x,k)$ parts:
       \begin{equation}
\wt \psi(x,k) =\psi_{\rm osc}  (x,k)  +\psi_{\rm dec} (x,k)
\label{eq:LSt}\end{equation}
(from now on, to simplify notation, we remove tilde   from  various objects related to the operator $\wt B$ and so we can write \e{eq:LSt} as  \e{eq:LS}).
It is natural to expect that the asymptotics of the integral \e{eq:BH} is determined by $\psi_{\rm osc}  (x,k)$ which yields the integral with the oscillating function $e^{  -i\xi (x)\ln t}e^{  i kx}$.

  \medskip

{\bf 6.3}.
    First, we will   show  that the decaying
  term $\psi_{\rm dec} (x,k)$ does not contribute to   the asymptotics of the integral \e{eq:BH}. To that end, we need the following elementary result.
 
  \begin{lemma}\label{StPhNo}
Let 
\[
{\sf J} (N)=     \int_{-\infty}^\infty e^{  -i\xi (x) N + i\kappa x } R(x) dx
\]
% \label{eq:BHN}\end{equation}
 where  $\kappa\in {\Bbb C}\setminus [0,\infty)$ $($or $\kappa\in {\Bbb C}\setminus (-\infty, 0])$ and
   \begin{equation}
   | e^{    i\kappa x } R^{(l)}(x) | \leq C (1+ |x|)^{-a-l}, \q l=0,1, \q a\in (0,1).
 \label{eq:BG3}\end{equation}
 Then ${\sf J} (N)=O( | N|^{-a} )$
 as  $N\to +\infty$ $($or  as  $N\to-\infty)$.
  \end{lemma} 
  
    \begin{pf}
    Integrating by parts we see that
 \begin{equation}
{\sf J} (N)=  i   \int_{-\infty}^\infty e^{ - i \xi(x) N +i \kappa x   }
\Big( \frac{R' (x)} {\kappa -\xi'(x) N}  + N \frac{R (x)\xi''(x)} {(\kappa-\xi'(x) N)^2}   \Big) d x.
 \label{eq:BF3}\end{equation}
Since $\xi'(x) >0$, it follows from our assumptions on $\kappa$ and $\sgn N$ that
  \begin{equation}
| \kappa-\xi'(x) N|\geq c(1+\xi'(x)| N|)  .
 \label{eq:BF3z}\end{equation}
Therefore  integral \e{eq:BF3} over $|x|\leq 1$ is bounded by $C |N|^{-1}$. According to \e{eq:BH6}
 expression \e{eq:BF3z} is minorated by $1+ |x|^{-1} |N |$ for $|x|\geq 1$.
  Using also  estimates \e{eq:BH7} and \e{eq:BG3} we see that      integral   \e{eq:BF3}  over $|x|\geq 1$ is bounded by
\[
  C   \int_1^\infty  
 x^{-a}(x+|N|)^{-1}   d x \leq C_{1}| N|^{-a}.
 \]
 %\label{eq:BF2x}\end{equation}
This yields the same estimate for the integral ${\sf J} (N)$.
     \end{pf}

  % \begin{pf}   It suffices to apply Lemma~\ref{StPhNo} with $\kappa=k$ and $R=w$.             \end{pf}

Recall that $\psi_{\rm dec}(x,k)$ is given by formula \e{eq:LS3}  where the function $ w  $ satisfies bound \e{eq:WFDw}.

 \begin{lemma}\label{StPhNo1}
We have the estimate 
  \begin{equation} 
{\cal I}_{\rm dec} (N,k): =   \int_{-\infty}^\infty e^{  -i\xi (x) N  } \zeta (x)  \psi_{\rm dec}(x,k)dx =O( | N|^{-1/2} )
      \label{eq:BHN1}\end{equation}
 as  $|N |\to \infty$.
  \end{lemma} 
  
    \begin{pf}
    Let us, for example, consider one of the terms in the first sum in \e{eq:LS3} and use Lemma~\ref{StPhNo}  with $\kappa=\kappa_{j}$ and
    \[
     R(x)= \zeta (x) \int_{-\infty}^x e^{  -i\kappa_{j} y } w(y,k) dy.
     \]
  Since $\Im \kappa_{j} >0$, it follows from \e{eq:WFDw} and \e{eq:BHNd1} that condition \e{eq:BG3}  is satisfied.
  This is obvious for $l=0$ and we have to use    estimate \e{eq:xd} for $l=1$.  Thus relation \e{eq:BHN1} is a consequence of  Lemma~\ref{StPhNo} for $a=1/2$.
          \end{pf}

       Next, we consider the contribution to  \e{eq:BH} of the oscillating  function  $\psi_{\rm osc}(x,k)$.  According to
        \e{eq:OSC2} for $n$ is even, it equals
           \begin{equation}
{\cal I}_{\rm osc} (N,k)=     \int_{-\infty}^\infty e^{  -i\xi (x) N  + ik x} \zeta (x) r_{+}(x,k)dx
+ \int_{-\infty}^\infty e^{  -i\xi (x) N  - ik x} \zeta (x) r_{-}(x,k)dx.
 \label{eq:ev}\end{equation}
 Here  $r_{+}(x,k)  $  and $r_{-}(x,k)  $ are given by  \e{eq:ev1}  and  \e{eq:ev3}, respectively.
  If $n$ is odd,  then according to
        \e{eq:OSC1}  formula  \e{eq:ev}\ is true with $r_{-}(x,k) =0$ and $r_{+}(x,k) = r(x,k)$    given by   \e{eq:OSC3}.  In all these formulas  $w(x,k) $ is defined by \e{eq:LS4} and hence according to \e{eq:WFDw}
   \[
r_{\pm} (x,k)=O (1), \q r_{\pm}'(x,k)=O (| x |^{-\rho}) \q {\rm as} \q |x|\to\infty.
\]
% \label{eq:theta}\end{equation}
 Using also \e{eq:BHNd1}, we see that
 \begin{equation}
     r_{\pm}(x,k) \zeta  (x) =O(|x|^{-1/2 }) , \q  ( r_{\pm}(x,k) \zeta  (x))'=O(|x|^{-3/2 }\ln\ln |x|) .
 \label{eq:BG2}\end{equation}

   Integrals \e{eq:ev} have stationary points
  defined by the equations 
  \[
N \xi' (x)=  k\q {\rm and}\q N \xi' (x)= - k,
 \]
 %\label{eq:STP}\end{equation}
 respectively. For $|N|$ is large enough,
 the first (the second) of these equations
  does not have solutions if $k N <0$ (if $k N > 0$), and they  have two solutions $\pm x_{N}$ if $k N >0$ for the first equation (or if $k N  < 0$ for the second equation). It easily follows from \e{eq:BH6} that
  \begin{equation}
x_{N}= |N| y_{N} \q {\rm where}\q y_{N}= \pi^{-1}  |k|^{-1}  n + a_{1}  |N|^{-1} + O (N^{-2}), \q |N|\to \infty.
 \label{eq:STP1}\end{equation}
 
 If there are no stationary points, then we have   the following result  which     is a particular case of Lemma~\ref{StPhNo}  corresponding to real $\kappa$.

  \begin{lemma}\label{StPhNo2}
If $|N|\to\infty$ and $\pm Nk<0$, then
 \[
  \int_{-\infty}^\infty e^{  -i\xi (x) N  \pm  ik x} \zeta  (x) r_{\pm}(x,k)dx= O( | N|^{-a} ), \q \forall a<1/2. 
 \]
 %\label{eq:BHNd}\end{equation}
  \end{lemma}

     The asymptotic behavior  of the  integrals \e{eq:ev} is of course determined by neighborhoods of the stationary points only.     Let us  check this statement.    
     Choose an even  function $\chi \in C_{0}^\infty ({\Bbb R})$ such that $\chi (y)=1$ if $|y|\leq \epsilon$ and  $\chi (y)=0$ if $|y|\geq 2\epsilon$  for a sufficiently small $\epsilon$ (for example, $\epsilon =(2\pi   |k|)^{-1}  n$)
 and set
\[
\wt{\chi} (y)= 1- \chi (y -y_{N})-\chi (y+y_{N}).
\]
The proof of the following result is very close to that of Lemma~\ref{StPhNo}.

   \begin{lemma}\label{StPo}
If $|N|\to\infty$ and $\pm Nk>0$, then
    \begin{equation}
     \int_{-\infty}^\infty e^{  -i\xi (x) N  \pm  ik x} \zeta  (x) r_{\pm}(x,k) \wt{\chi} (x/| N|) dx
= O(|N|^{-a}) , \q \forall a<1/2. 
 \label{eq:stp1}\end{equation}
  \end{lemma} 

\begin{pf}
   Integrating by parts we see that integral \e{eq:stp1} equals integral \e{eq:BF3} where $\kappa= \pm k$ and
   \[
   R(x)= \zeta  (x) r_{\pm}(x,k) \wt{\chi} (x/|N|).
   \]
    Since $\wt{\chi} (y) =0$ for $y\in (y_{N}-\epsilon, y_{N}+\epsilon)$ and $y\in (-y_{N}-\epsilon, -y_{N}+\epsilon)$, in the integral \e{eq:BF3} we now have either $|x|\leq |N| (y_{N}-\epsilon)$ or $|x|\geq |N| (y_{N} +\epsilon)$. In the first and second cases $ |\pm k-\xi'(x)| N|$ is minorated by $c (1+|x|)^{-1} |N|$ and a positive constant $c$, respectively. Therefore   using   relations \e{eq:BH6} and \e{eq:BH7} and  estimates \e{eq:BG2}, we conclude the proof of \e{eq:stp1} quite similarly to   Lemma~\ref{StPhNo}.
    \end{pf}

 %\begin{equation}
%   i  \int_{-\infty}^\infty e^{  -i\xi (x) N  \pm  ik x}  \Big( \frac{\big(\zeta  (x) r_{\pm}(x,k) \wt{\chi} (x/|N|)\big)'} {\pm k-\xi'(x) N}  + N \frac{ \zeta (x) r_{\pm}(x,k) \wt{\chi} (x/ |N|)\xi''(x)} {(\pm k-\xi'(x) N)^2}   \Big) d x.
% \label{eq:BFz}\end{equation}
% Below we use relations \e{eq:BH6} and \e{eq:BH7} and  estimates \e{eq:BG2}.
%    Since $\wt{\chi} (y) =0$ for $y\in (y_{N}-\epsilon, y_{N}+\epsilon)$ and $y\in (-y_{N}-\epsilon, -y_{N}+\epsilon)$, in the integral \e{eq:BF3} we have either $|x|\leq |N| (y_{N}-\epsilon)$ or $|x|\geq |N| (y_{N} +\epsilon)$.   In the first case, for some $c>0$,   the estimate 
% \[|k|-\xi'(x)| N| \geq c \zeta'(x)| N| \geq c_{1} (1+|x|)^{-1/2}| N|\]
%holds. Therefore the first of these integrals is estimated by
%\[| N|^{-1}\int_{0}^{ |N| (y_{N}-\epsilon)}(1+|x|)^{-1/2} dx =C | N|^{-1/2}.\]
%   In the second case we have   the estimate 
% \[ \zeta'(x)| N| -|k|\geq c >0.\]
%Therefore the second integral is estimated by
% \[ \int_{ |N| (y_{N}+\epsilon)}^\infty (1+|x|)^{-3/2} dx =C | N|^{-1/2}.\]
%This concludes the proof of \e{eq:stp1}.

     \medskip

{\bf 6.4}.
Finally, we consider neighborhoods of the stationary points $y_{N}$ defined by \e{eq:STP1}.

 \begin{lemma}\label{SxPH}
Put 
   \begin{equation}
 \phi(y,k; N)= -\xi (|N| y) + | k| y .
 \label{eq:BH2}\end{equation} 
 Then, as $|N| \to \infty$ and $\pm N k>0$, the integral \e{eq:ev} satisfies the asymptotic relation
   \begin{multline}
{\cal I}_{\rm osc} (N,k)=  (2 n)^{1/2}   | k|^{-1}  | N|^{1/2}  \sum_{\tau=``\pm"} e^{ \pi i \tau  \sgn N  /4 }    e^{i N \phi  (\tau y_{N},k; N)}   
\\
\times \zeta (\tau  | N| y_{N} ) r_\pm( \tau  | N| y_{N}   ,k)     + O(|N|^{-a}  ), \q \forall a<1/2.
 \label{eq:BF}\end{multline}
 \end{lemma} 
 
  \begin{pf}
  It follows from Lemmas~\ref{StPhNo2} and \ref{StPo}  that
    \begin{equation}
   {\cal I}_{\rm osc} (N,k)=      {\cal I}_{\pm} (N,k) + O (|N|^{-a}), \q \pm N k>0,
 \label{eq:pm}\end{equation}
where
    \[
   {\cal I}_{\pm} (N,k)=    \int_{-\infty}^\infty e^{  -i\xi (x) N  \pm  ik x} \zeta (x) r_{\pm}(x,k) \big(\chi (x/| N|- y_{N})
   + \chi (x/| N|+ y_{N}) \big) dx.
 \]
 %\label{eq:stp2}\end{equation}
 Let us find the asymptotics of this integral.
 Making the change of variables $x=|N|y$, we see that 
  \begin{equation}
   {\cal I}_{\pm} (N,k)=    |N | \sum_{\tau=``\pm"} {\cal J}_{\pm}^{(\tau)} (N,k)
 \label{eq:IJ}\end{equation}
   where 
                \[
{\cal J}_{\pm}^{(\tau)} (N,k) =     \int_{-\infty}^\infty e^{  iN \phi(y,k; N) } \zeta  ( | N| y) r_\pm ( | N| y,k)  \chi (y -\tau  y_{N})  dy .
 \]
 %\label{eq:AA}\end{equation}  
 % We consider two integrals in the right-hand side of \e{eq:AA} separately.
After the change of variables
   $x=y - \tau y_{N}$, we find that 
       \begin{equation}
 {\cal J}_{\pm}^{(\tau)} (N,k) = \int_{-\infty}^\infty e^{  iN \phi(x + \tau y_{N} ,k; N) } \zeta (| N|(x + \tau y_{N})  )r_\pm ( | N|(x+ \tau y_{N}) ,k)   \chi (x )    dx .
 \label{eq:BFx}\end{equation}
 
  Let us apply to this integral Lemma~\ref{StPh} with the functions
 \begin{equation}
 \omega  (x,k; N)= \phi ( x +\tau  y_{N} , k) 
  \label{eq:BFxh}\end{equation}
 and
  \[
 g  (x,k; N)= \zeta (| N| (x +\tau y_{N}) ) r_\pm( | N|(x +\tau y_{N})   ,k)   \chi (x )
 \]
  depending on the parameter $N$ (as well as on $\tau$ and $k$). 
Since we integrate in \e{eq:BFx} over   $\Bbb R$ (instead of  $\Bbb R_{+}$ in \e{eq:BHs}), we have to multiply   the right-hand side of \e{eq:BHs1} by the additional factor $2$. Let us check  the assumptions of Lemma~\ref{StPh}.
 Differentiating formula \e{eq:BFxh} $p$ times in $x$ and using \e{eq:BH2}, we see that 
  \[
 \omega^{(p)} (x,k; N)= - |N|^p \xi^{(p)}(|N |(x +\tau y_{N})),\q p\geq 2. 
\]
% \label{eq:AA1}\end{equation}
In particular, according to \e{eq:BH7}, \e{eq:STP1} we have
    \begin{equation}
 \omega ''    (0,k; N)= -N^2 \xi'' (\tau | N | y_{N})= \tau \pi^{-1} n y_{N}^{-2}+ O(|N|^{-1})=
 \tau \pi  k^2 n^{-1}+ O(|N|^{-1}).
 \label{eq:BF1}\end{equation}
  Moreover, all derivatives of $ \omega  (x, k;N)$ are uniformly in $N$ bounded for $|x|\leq 2 \epsilon$ because  $\xi^{(p)} (x)=O (|x|^{-p})$ as $|x|\to\infty$  (see \e{eq:BH7}).
We have
    \begin{equation}
 g  (0,k; N)=  \zeta (  \tau | N|  y_{N} ) r_\pm( \tau | N|  y_{N}  ,k) ,  
 \label{eq:Bg}\end{equation}
 and it follows from estimates    \e{eq:BG2}  that the functions $g^{(p)} (x,k; N)=O( |N|^{-a})$ for $p=0, 1$ and arbitrary  $\varepsilon>0$. Therefore the remainder $ {\cal R}_{\pm}^{(\tau)} (N,k)$ in formula  \e{eq:BHs1} for $ {\cal J}_{\pm}^{(\tau)} (N,k) $ is bounded by $C |N|^{-3/2+\varepsilon}  $.  Now \e{eq:BHs1} and equalities \e{eq:BF1}, \e{eq:Bg} imply that
    \begin{multline*}
 {\cal J}_{\pm}^{(\tau)} (N,k)= (2 n)^{1/2}   | k|^{-1}  | N|^{-1/2} e^{\pi i  \tau  \sgn N  /4 }    e^{i N \phi  (\tau y_{N},k; N)}   
 \\
\times  \zeta  ( \tau | N| y_{N}  ) r_\pm( \tau | N| y_{N}  ,k)    + O(|N|^{-1-a}  ) .
% \label{eq:BGg}
\end{multline*}
 
 Substituting this expression into \e{eq:IJ} and taking into account relation \e{eq:pm}, we conclude the proof of \e{eq:BF}.
     \end{pf}

It remains to study the right-hand side of \e{eq:BF}.

   \begin{lemma}\label{StPH}
       Let  the function $\gamma (N,k)$ and the  integral  ${\cal I}  (N,k)$     be given by equalities  \e{eq:BFGT} and  \e{eq:BH}, respectively. 
Then as $|N| \to \infty$, $\pm Nk>0$,    the asymptotic formula
   \begin{equation}
{\cal I}  (N,k)= \sqrt{\frac{2n}{ |k|}} \big( e^{i \gamma(N,k)} r_{\pm} (+\infty,k) + e^{-i \gamma(N,k)} r_{\pm} (-\infty,k)\big)  (1+ O(|N|^{-a}))
 \label{eq:As}\end{equation}
holds with an arbitrary $a<1/2$.
 \end{lemma} 
 
    \begin{pf}
 First, we note that in view of Lemma~\ref{StPhNo1}  it suffices to prove \e{eq:As} for the integral ${\cal I}_{\rm osc}  (N,k)$.
Therefore we only have to 
  find the asymptotic behavior of the right-hand side of \e{eq:BF} as $|N|\to\infty$. We start with the phase  
   $ \phi ( \pm y_{N},k; N)$ given by relation \e{eq:BH2}.   
    Using   asymptotic formulas \e{eq:BH5} and  \e{eq:STP1},  we see that  
  \[
 \phi( y_{N},k; N)= -\phi(- y_{N},k; N)= \gamma_{0} (| N /k|) + a_{1}|k| |N|^{-1} + O(N^{-2}) .
\]
% \label{eq:STP2}\end{equation}
 Next, we consider the function $  \zeta (\pm |N| y_{N}  ) $ defined by  formula \e{eq:BH1}. 
It follows from 
  \e{eq:han4q} and  \e{eq:STP1} that its phase
  \[
  \varrho ( |N| y_{N}  ) =-  \varrho (- |N| y_{N}  ) =  \gamma_1 (| N /k|)   + O(N^{-1}).
  \]
  The asymptotics of its modulus is a direct consequence of  \e{eq:BH6}:
        \[
        \xi' ( \pm |N| y_{N}  ) = |k/N|^{1/2}     \big( 1+ O (N^{-1})\big) .
   \]

 Let us also use   that, by definition  \e{eq:ev1}, there exists  finite limits $r_\pm (+\infty,k)$ and $r_\pm (-\infty,k)$.
 Substituting the results obtained into the right-hand side of \e{eq:BF}, we conclude the proof of \e{eq:As}.
 \end{pf}
 
  Now we are in a position to conclude the {\it proof }  of Theorem~\ref{SpThHb} and hence of Theorem~\ref{SpThHa}. 
  We proceed from relation \e{eq:BHC} and use that $\omega (t,k)= \gamma (\ln t,k)$.
 So we only have to calculate the limits $r_{\pm} (+\infty,k) $ and $r_{\pm} (-\infty,k) $ in the right-hand side of \e{eq:As}.  
  
   If  $n$ is odd, then  equality  \e{eq:OSC1} ensures that $r_{+} (x,k)= r (x,k)$ and $r_{-} (x,k)= 0$.  According to \e{eq:OSC3}  we have $ r (-\infty,k) =1$ and according to \e{eq:DES}    $r  (+\infty,k)=s (\lambda)$. 
   This yields assertion (i)    of Theorem~\ref{SpThHb}.

Let   $n$ be even.  Recall that, by the definition of the scattering matrix in subs.~2.4,  if $k>0$, then  $r_{+} (+\infty,k)=s_{11}(\lambda)$ according to \e{eq:SMX} and $ r_{+} (-\infty,k) = 1$ according to \e{eq:ev1}. Therefore  formula \e{eq:As}   implies the  asymptotic relation \e{eq:S1H}.
If $k <0$, then  $r_{+} (+\infty,k)=1$ according to \e{eq:ev1} and $ r_{+} (-\infty,k) = s_{22}(\lambda)$ according to \e{eq:SMX}. This leads to the   asymptotic relation \e{eq:S2H} and concludes the proof of assertion (ii)    of Theorem~\ref{SpThHb}. $\qed$

   \section{Large times  evolution}
   
   Here we discuss the asymptotic behavior  of $\exp(-iHT)u$ as $T\to \pm\infty$ for $u$  in the absolutely continuous subspace ${\cal H}^{(ac)}_{H} $ of the generalized Carleman operator $H$. The construction below is rather similar to the proof of Theorem~\ref{SpThHb}, and so we omit some technical details. We have to distinguish the cases of even and odd $n$.

    \medskip

{\bf 7.1}.
It follows  from relations  \e{eq:V} and  \e{eq:MM1}  that
\[
\exp(-iHT)u= F^* \mathrm{L}\exp(-i\ B T)\mathrm{L}^* F u
\]
where the asymptotics of $\exp(-i\ B T)$ is given by Theorem~\ref{qqs}. According to Proposition~\ref{BL} and equality   \e{eq:gauge} the functions 
$\Theta$ and $ \Theta_{1}$, $\Theta_{2}  $ in formulas \e{eq:StPhB} and \e{eq:StPhBod}  now equal $\beta(x)=-q_{n-1} \xi(x)/n$.
{\it In calculations below, we neglect the terms whose norms in $L^{2}({\Bbb R})$ tend to zero as }
$T\to \pm\infty$. 

Suppose first that $n$ is even.   Using  the definitions  \e{eq:V1} and  \e{eq:MAID1} of the operators $\mathrm{L}$ and $F$ and the asymptotic formula  \e{eq:StPhB}, we see that 
 \begin{multline}
(\exp(-iHT)u)(t)=  | 2\pi  n T |^{-1/2}  e^{\mp \pi i/4}
\\
\times \int_{-\infty}^\infty  t^{-1/2-i\xi} e^{i\eta(\xi)-iq_{n-1}  \xi /n} e^{i  \Xi_{0}(x(\xi), T)}x'(\xi)^{1/2} f_{\pm}(\frac{x({\xi})}{nT} )d\xi
 \label{eq:AST}\end{multline}
 where the phase $ \Xi_{0}(x, T)$ is given by formula \e{eq:Xi} and
  \begin{equation}
 f_{\pm}= YW_{\pm} (B_{0}, B ;\mathcal{J}^*)\mathrm{L}^* F u
 \label{eq:ASTW}\end{equation}
 with $Y$ defined by \e{eq:StPh1}. Our goal is to find the asymptotics as $T\to\pm \infty$ of integral \e{eq:AST} in $L^2 ({\Bbb R})$;  so we may suppose that $f_{\pm} $ belongs to its dense subset $ C_{0}^\infty ({\Bbb R}\setminus \{0\})$. 
 Let us make   the changes of variables $\xi\mapsto x=x(\xi)$ and then $x\mapsto y= x /(nT)$. Using notation \e{eq:han4}, \e{eq:BH1}, we rewrite \e{eq:AST} as
  \begin{equation}
(\exp(-iHT)u)(t)=  \sqrt{\frac {n  |T |}{ 2\pi t}}  e^{\mp \pi i/4}
 \int_{-\infty}^\infty   e^{i  (-\xi (nTy) \ln t +(n-1) |y|^{\frac{n}{n-1}} T )}\zeta(nTy)  f_{\pm}( y )d y.
 \label{eq:AST1}\end{equation}
 
 We can now apply the stationary phase method (see, e.g., Appendix~A where the role of $T$ is played by $N$) with the phase function
  \begin{equation}
\omega (y, t, T)=-\xi (nTy) T^{-1} \ln t +(n-1) | y|^{\frac{n}{n-1}} 
 \label{eq:AST3}\end{equation}
 depending on the parameters $t$ and $T$. The stationary points  $y= y(t,T)$ of this function are determined by the equation $\omega'_{y} (y, t, T)=0$, i.e.,
   \begin{equation}
   \xi' (nTy) \ln t =|y|^{\frac{1}{n-1}} \sgn y .
 \label{eq:AST2}\end{equation}
Recall that the function $\xi'(x)$ is even,  $\xi'(x)> 0$, and it has asymptotics  \e{eq:BH6} as $x\to+\infty$. Therefore for sufficiently large $| T|$, equation \e{eq:AST2} has the unique solution and
   \begin{equation}
  y (t,T)= \sgn (\ln t) \Big|\frac{\ln t} {\pi T}\Big|^{\frac{n-1}{n}}  \big(1+ O (|T|^{-1})\big).
 \label{eq:AST4}\end{equation}
 Differentiating  \e{eq:AST3} twice and using asymptotics \e{eq:BH7}, we see that
 \begin{align*}
\omega''_{yy} (y, t, T)=&-n^2 \xi '' (nTy) T  \ln t +n(n-1)^{-1} | y|^{-\frac{n-2}{n-1}}
\\
=&n     (\ln t/(\pi T)) y^{-2}  \sgn (Ty) \big(1+ O (|Ty|^{-1}) \big)+n(n-1)^{-1} | y|^{-\frac{n-2}{n-1}}  
 %\label{eq:AST5}
 \end{align*}
 as $|Ty|\to\infty$. Substituting here expression \e{eq:AST4}, we find that
 \begin{equation}
\omega''_{yy} ( y (t,T), t, T)
=n^2 (n-1)^{-1} \Big|\frac{\ln t} {\pi T}\Big|^{-\frac{n-2}{n}}  \big(1+ O (|T|^{-1})\big)
 \label{eq:AST6}\end{equation}
 provided $T^{-1}  \ln t$ belongs to compact subsets of ${\Bbb R}\setminus\{0\}$.
 
 Note also that according to  \e{eq:BH1} and \e{eq:BH6} we have
 \begin{equation}
 \zeta (nT y(t,T))= e^{i \varrho (nT  y (t,T)) }  (\pi |T|)^{-1/2}\, \Big| \frac{\ln t} {\pi T} \Big|^{- (n-1)/(2n)} \big(1+ O (|T|^{-1})\big).
 \label{eq:AST6a}\end{equation}

Thus formula  \e{eq:BHs1} yields the following assertion.
 
  \begin{theorem}\label{HTU}
Let $H$ be the Hankel operator \e{eq:H1bis} with  kernel      \e{eq:LOGbis} where $n$ is even. Let $y(t,T)$ be the $($unique$)$ solution of equation  \e{eq:AST2}. Set 
 \begin{equation}
\varPhi ( t, T)=-\xi (nT  y (t,T) )  \ln t +(n-1) |  y (t,T) |^{\frac{n}{n-1}} T
 + \varrho (nT  y (t,T))  .
 \label{eq:AST7}\end{equation}
Then for all  $u\in{\cal H}^{(ac)}_{H} $   the asymptotic relation holds
   \begin{multline*}
   (\exp(-iHT)u)(t)=e^{i \varPhi ( t, T)} \sqrt{ \frac{n-1}{\pi n t |T|}} \, \Big| \frac{\ln t} {\pi T} \Big|^{-1/(2n)}
   \\
   \times f_{\pm}\big(\sgn (\ln t)  \Big| \frac{\ln t} {\pi T} \Big|^{(n-1)/n} \big)
   +\epsilon_{\pm}(t,T)
% \label{eq:AST8}
 \end{multline*}
 where    $f_{\pm}$ is defined by formula \e{eq:ASTW} and $\| \epsilon_{\pm}(T)\|\to 0$ as $T\to\pm\infty$.
   \end{theorem}
 
 \medskip

{\bf 7.2}.
 Let us now briefly consider the case of odd $n >1$. We proceed again from Theorem~\ref{qqs} but use formula \e{eq:StPhBod}. Therefore instead of  \e{eq:AST1} we have the representation
   \begin{equation}
(\exp(-iHT)u)(t)=    (  n  |T |/ (2\pi t))^{1/2}\sum_{j=1,2} e^{\mp (-1)^{j }\pi i/4}
 \int_{0}^\infty   e^{i  \omega_{j} (y,t,T) T }\zeta(nTy)  f_{\pm}( (-1)^{j} y )d y
 \label{eq:AST1-}\end{equation}
 where  $T\to\pm\infty$,
  \[
\omega_1 (y, t, T)=-\xi (nTy) T^{-1} \ln t - (n-1)  y^{\frac{n}{n-1}} ,\q y>0,
 \]
 % \label{eq:AST3-}\end{equation}
 and $\omega_2 = \omega $ is given by formula \e{eq:AST3} with $y>0$. The stationary points of the integral in \e{eq:AST1-} for $j=1$ are determined by the equation $\partial\omega_1 (y,t,T) /\partial y =0$, that is,
   \begin{equation}
    \xi' (nTy) \ln t =-y^{\frac{1}{n-1}} .
  \label{eq:stop1}\end{equation}
   Let $| T|$ be sufficiently large. Then this  equation  has no solutions for $t>1$, and it has the unique solution $ y_{1} (t,T)$ for $t<1$. In view of \e{eq:BH6} this solution satisfies the asymptotic relation
   \[
  y_{1} (t,T)= \Big|\frac{\ln t} {\pi T}\Big|^{\frac{n-1}{n}}  \big(1+ O (|T|^{-1})\big).
  \]
% \label{eq:stop}\end{equation}
 Instead of \e{eq:AST6} we now have
 \[
\partial^2\omega_1 ( y_{1} (t,T),t,T) /\partial y^2   
=- n^2 (n-1)^{-1} \Big|\frac{\ln t} {\pi T}\Big|^{-\frac{n-2}{n}}  \big(1+ O (|T|^{-1})\big).
 \]
  According to formula   \e{eq:AST3} for the function $\omega_2 ( y,t,T)$,  the stationary points of the integral in \e{eq:AST1-} for $j=2$ are determined by the equation   \e{eq:AST2} where  $y>0$. For   $| T|$   sufficiently large, this  equation  has no solutions for $t<1$, and it has the unique solution $ y_2 (t,T)$ for $t >1$. This solution satisfies the asymptotic relation \e{eq:AST4} (where $ \sgn (\ln t)=1$). The second  derivative of $\omega_2 ( y,t,T)$ satisfies relation   \e{eq:AST6}.

 Now using  formula  \e{eq:BHs1} and taking into account asymptotic relation \e{eq:AST6a}  we obtain the following assertion.
 
  \begin{theorem}\label{HTUodd}
Let $H$ be the Hankel operator \e{eq:H1bis} with  kernel      \e{eq:LOGbis} where $n>1$ is odd. Let $y_{1}(t,T)$ for $t<1$ and  $y_2(t,T)=y(t,T)$ for $t>1$ be the $($unique$)$ solutions of equations  \e{eq:stop1} and \e{eq:AST2}, respectively. Put 
 \[
\varPhi_{1} ( t, T)=-\xi (nT  y_{1} (t,T) )  \ln t -(n-1)   y_{1} (t,T) ^{\frac{n}{n-1}} T 
 + \varrho (nT  y_{1} (t,T)) ,  \q t<1, 
 \]
 and let $\varPhi_{2}( t, T)=\varPhi  ( t, T)$ be defined for $t>1$ by formula \e{eq:AST7}.
Then for all  $u\in{\cal H}^{(ac)}_{H} $   the asymptotic relation holds
   \begin{multline*}
   (\exp(-iHT)u)(t)= \sqrt{ \frac{n-1}{\pi n t |T|}} \, \Big| \frac{\ln t} {\pi T} \Big|^{-1/(2n)}
   \\
   \times \sum_{j=1,2} e^{i \varPhi_{j} ( t, T)}  \chi_{j} (t) f_{\pm}\big( (-1)^j \Big| \frac{\ln t} {\pi T} \Big|^{(n-1)/n} \big)
   +\epsilon_{\pm}(t,T)
% \label{eq:AST8odd}
 \end{multline*}
 where $\chi_1 (t) $, $\chi_2 (t) $ are the characteristic functions of the intervals $(0,1)$ and $(1,\infty)$, respectively,  $f_{\pm}$ is defined by formula \e{eq:ASTW} and $\| \epsilon_{\pm}(T)\|\to 0$ as $T\to\pm\infty$.
   \end{theorem}
   
    Theorems~\ref{HTU} and \ref{HTUodd} show that  as $T\to \pm\infty$  the function $(\exp(-iHT)u) (t)$ ``lives" in
    exponentially small neighborhoods of the singular points $t=0$ and $t=\infty$, that is, in the region where $t \sim e^{\pi |T|}$ or $t \sim e^{-\pi |T|}$. This is quite different from differential operators \e{eq:VB}  when according to Theorem~\ref{qqs} we typically have  $|x|\sim n|T|$ for the evolution operator $(\exp(-iBT)f) (x)$.

%{\bf 7.3}.
% Finally, we compare Theorems~\ref{HTU} and \ref{HTUodd}
% with Theorem~\ref{qqs} where the asymptotic behavior as $T\to \pm\infty$ of the Schr\"odinger evolution $\exp(-iBT)f$ was described. If $n$ is even, then according to \e{eq:StPhB} the function $(\exp(-iBT)f) (x)$ ``lives" in the region where $x\sim n|T|$. If $n$ is odd, then according to \e{eq:StPhBod} this result remains true but additionally $(\exp(-iBT)f) (x)$ quits the half-axis ${\Bbb R}_{\mp}$.
 
 %Theorems~\ref{HTU} and \ref{HTUodd} show that the behavior of the function $(\exp(-iHT)u) (t)$ is of the same nature.
%This function   ``lives" in the region where $t \sim e^{\pi |T|}$ or $t \sim e^{-\pi |T|}$.
  
\appendix
%%%%%%%%%%%%%%%%%%%%%%%%%%%%%%%%%%%%%%%%%%%
%%%%%%%%%%%%%%%%%%%%%%%%%%%%%%%%%%%%%%%%%%%
  %%%%%%%%%%%%%%%%%%%%%%%%%%%%%%%%%%%%%%%%%%%
%%%%%%%%%%%%%%%%%%%%%%%%%%%%%%%%%%%%%%%%%%%
\section{Stationary phase method}

 The statement below is quite standard. We give its short proof because apparently explicit estimates of the remainders are not available in the literature.  
  
     \begin{lemma}\label{StPh}
Let 
\begin{equation} 
 {\cal J}  (N)=       \int_0^\infty e^{  iN \omega (y) }
g (y) d y
 \label{eq:BHs}\end{equation}
 where  $\omega =\bar{\omega}\in C^\infty ({\Bbb R})$, $\omega'(0)=0$ and $g\in C_{0}^\infty ({\Bbb R})$. Suppose that $g(y)=0$ if $y \geq a$ for some $a>0$ and $\omega''(y) \neq 0$ for $y\in [ 0,a ]$. Put
 \[
 g_{n}=\max_{y\in [0,a]}| g^{(n)} (y)|, \q  \omega_{n}=\max_{y\in [0,a]}| \omega^{(n)} (y)|, \q 
\kappa =\min_{y\in [0,a]} |\omega'' (y)|.
 \]
 Then
 \begin{equation}
{\cal J} (N)=   2^{-1} e^{   i\tau \pi/4 + i N \omega (0)}    (2\pi )^{1/2}   | \omega''(0) N|^{-1/2}  g(0)  +
{\cal R} (N), \q | N|\to \infty,
 \label{eq:BHs1}\end{equation}
 where $\tau=\sgn(\omega''(0) N)$ and
  \begin{equation} 
|{\cal R} (N)|\leq C  (\kappa^{-7/2} \omega_{2}^{3/2}\omega_{3}  g_{0}
+ \kappa^{-2}   \omega_2  g_{1})
| N |^{-1 } \big(1+ \big|\ln | \omega_{0}  N | \big|\big)
 \label{eq:BHs2}\end{equation}
 with some absolute constant $C$.
  \end{lemma} 

\begin{pf}
%Suppose for definiteness that $N\to+\infty$.
Without loss of generality, we may suppose that $\omega(0)=0$ and $\omega''(0) >0$.
Let us make the change of variables $x=\omega(y)$ in \e{eq:BHs}. Then 
\begin{equation}
{\cal J} (N)=       \int_0^\infty e^{  iN x }
{\bf g}(x)  x^{-1/2}d x
 \label{eq:BHs3}\end{equation}
 where 
  \begin{equation}
 {\bf g}(x) =\sigma (\eta(x)) g (\eta(x)) ,  
 \label{eq:BHs4}\end{equation}
 \begin{equation}
  \sigma(y)= \omega(y)^{1/2} \omega'(y)^{-1}
 \label{eq:BHs8}\end{equation}
  and $\eta (x)$ is the function inverse   to $\omega(y)$. Integrating in \e{eq:BHs3} by parts, we see that
  \begin{equation}
{\cal J} (N)=      {\bf g}(0) \int_{0}^\infty e^{  iN y }
  y^{-1/2}d y +
 \int_0^\infty \big(  \int_x^\infty e^{  iN y } y^{-1/2}dy\big)
{\bf g}'(x) d x.
 \label{eq:BHs5}\end{equation}
 According to \e{eq:BHs4} we have
 \[
 {\bf g}(0)= g(0)   \sqrt{\lim_{y\to 0}\omega(y) \omega'(y)^{-2}}= g(0)    (2 \omega''(0) )^{-1/2},
 \] 
 and hence the first terms in the right-hand sides of \e{eq:BHs1} and \e{eq:BHs5} coincide. 
 
 Let ${\cal R} (N)$ be the second term in the right-hand side of \e{eq:BHs5}. Since 
 \[
 \big|\int_x^\infty e^{  iN y } y^{-1/2}dy\big| \leq C  |N|^{-1/2} (1+ |N |x)^{-1/2}  ,
 \]
  we get the estimate
   \begin{equation}
|{\cal R} (N)| \leq C |N|^{-1/2 }
 \int_0^{\omega_{0}}  |{\bf g}'(x)| (1+ |N |x)^{-1/2}dx .
 \label{eq:BHs6}\end{equation}
 
 According to \e{eq:BHs8} we have $ \eta' (x) = \omega' (y)^{-1}=   \sigma(y) x^{-1/2}$. Therefore
 differentiating \e{eq:BHs4}, we see that
  \begin{equation}
 {\bf g}'(x) =(\sigma' (y) g (y) +\sigma (y) g' (y) )\sigma(y) x^{-1/2} , \q  y=\eta (x).
 \label{eq:BHs7}\end{equation}
  Observe that
 \begin{equation}
|\omega (y)|\leq \omega_{2} y^2/2,  \q \omega' (y)=\int_{0}^y \omega'' (x)d x\geq \kappa y  ,
\q   \omega(y)=\int_{0}^y \omega' (x)dx\geq \kappa y^2 /2,
 \label{eq:BHr}\end{equation}
 and hence, by definition \e{eq:BHs8},
  \begin{equation}
| \sigma (y)|\leq (\omega_{2}/2)^{1/2}  \kappa^{-1}  .
 \label{eq:BHr1}\end{equation}
 Note that
    \begin{equation}
\sigma' (y)=\frac{\omega'(y)^2 -2 \omega(y) \omega''(y)} {2 \omega(y)^{1/2}\omega'(y)^2}.
 \label{eq:BHs9}\end{equation}
 It  follows from  the first estimate \e{eq:BHr}  that $((\omega')^2 -2 \omega\omega'')'= -2 \omega\omega'''$  is bounded  by $\omega_{2}\omega_{3} y^2$,  and so the numerator in \e{eq:BHs9} is bounded  by $\omega_{2}\omega_{3} y^3 /3$.
 A lower bound by $2^{ 1/2}   \kappa^{5/2} y^3$  on the denominator in \e{eq:BHs9} follows from  the second and third estimates \e{eq:BHr} whence
     \begin{equation}
| \sigma' (y)|\leq  2^{- 1/2}  3^{-1} \kappa^{-5/2}  \omega_{2}\omega_{3} .
 \label{eq:BHs9c}\end{equation}
 Let us come back to function \e{eq:BHs7}.
Combining   \e{eq:BHr1} and  \e{eq:BHs9c}, we get the bound
  \[ 
|{\bf g}' (x)|\leq C  (\kappa^{-7/2}\omega_{2}^{3/2}\omega_{3}  g_{0}
+ \kappa^{-2}  \omega_2  g_{1})x^{-1/2} .
 \]
 In view of \e{eq:BHs6}, this yields \e{eq:BHs2}.
\end{pf}

%%%%%%%%%%%%%%%%%%%%%%%%%%%%%%%%%%%%%%%%
%%%%%%%%%%%%%%%%%%%%%%%%%%%%%%%%%%%%%%%%

 \end{document}